\documentclass[11pt,a4paper]{article}%
\usepackage{amscd}
\usepackage{amsmath}
\usepackage{graphicx}
\usepackage{amsfonts}
\usepackage{amssymb}
\usepackage{xcolor}
\usepackage{amsbsy,amsthm}
\usepackage{pifont}
\usepackage{fancyhdr,ifthen}
\usepackage[colorlinks=true, bookmarksnumbered=true, bookmarksopen=true,
bookmarksopenlevel=3, pdfstartview=FitH, linkcolor=cyan, pdfmenubar=true,
pdftoolbar=true, bookmarks=true,citecolor=cyan, urlcolor=magenta,
filecolor=cyan,menucolor=black,plainpages=false,pdfpagelabels]{hyperref}%
\newtheorem{theorem}{Theorem}[section]
\newtheorem{lemma}[theorem]{Lemma}

\newtheorem{proposition}[theorem]{Proposition}
\newtheorem{definition}[theorem]{Definition}
\newtheorem{remark}[theorem]{Remark}

\numberwithin{equation}{section}

\begin{document}

\title{Global solutions for a 2D chemotaxis-fluid system with large measures as
initial density and vorticity \vspace{0.5cm} }
\author{\textbf{Lucas C. F. Ferreira}\thanks{LCFF was supported by CNPq (grant: 308799/2019-4 and grant: 312484/2023-2), Brazil.}\\{\small Universidade Estadual de Campinas, IMECC,}\\{\small CEP 13083-859, Campinas-SP, Brazil.}\\{\small \texttt{email:\ lcff@ime.unicamp.br}}\vspace{0.8cm} \\\textbf{Daniel P. A. Lima}\thanks{DPAL was supported by CAPES (Finance Code 001) and FAPESP (grant: 2019/02733-1), Brazil.}\\{\small Universidade Estadual de Santa Cruz, DCEX,}\\{\small CEP 45662-900, Ilheus-BA, Brazil.}\\{\small \texttt{email:\ dpalima@uesc.br}}}
\date{}
\maketitle

\begin{abstract}
We consider a chemotaxis-fluid system in the whole plane $\mathbb{R}^{2}$
which describes the motion of bacteria suspended in a Navier-Stokes fluid and
attracted by a chemical (oxygen). Employing the vorticity formulation for the
fluid equations and considering the prototypical case, we obtain local and
global solutions with large (Radon) measures as initial data for the bacterial
density and vorticity. The gravitational/centrifugal potential is taken with
finite $L^{2}$-gradient that can be large. The uniqueness property is also
discussed. For the global result, we need to assume a smallness condition only
on the $L^{\infty}$-norm of the initial oxygen concentration. In comparison
with previous works, our results provide a new class for the initial density
and vorticity, as well as for the potential, covering particularly singular
measures such as Dirac delta, measure concentrated on smooth curves (filaments
and rings), among others. For that, we approach the system via critical
functional spaces, Kato-type norms, and suitable $L^{p}$-estimates uniformly
in time.



\end{abstract}

\pagestyle{fancy}
\fancyhf{}
\renewcommand{\headrulewidth}{0pt}
\chead{\ifthenelse{\isodd{\value{page}}}{L.C.F. Ferreira and D.P.A. Lima}{Chemotaxis-Navier-Stokes system with initial measures}}
\rhead{\thepage}

\section{Introduction}

In this paper we are concerned with the chemotaxis-Navier-Stokes model
\begin{equation}%
\begin{cases}
\displaystyle\partial_{t}n+u\cdot\nabla n=\Delta n-\nabla\cdot(n\chi(c)\nabla
c), & \qquad x\in\mathbb{R}^{2}\ ,\ t>0,\ \\[3mm]%
\partial_{t}c+u\cdot\nabla c=\Delta c-nf(c), & \qquad x\in\mathbb{R}%
^{2}\ ,\ t>0,\ \\[3mm]%
\partial_{t}u+u\cdot\nabla u=\Delta u+\nabla P+n\nabla\phi, & \qquad
x\in\mathbb{R}^{2}\ ,\ t>0,\ \\[3mm]%
\nabla\cdot u=0, & \qquad x\in\mathbb{R}^{2}\ ,\ t>0,\ \\[3mm]%
n(x,0)=n_{0}(x),\ \ c(x,0)=c_{0}(x),\ \ u(x,0)=u_{0}(x), & \qquad
x\in\mathbb{R}^{2},
\end{cases}
\label{CNSC}%
\end{equation}
where $n\geq0$ is the unknown bacterial density, $c\geq0$ is the oxygen
concentration, $u$ is the fluid velocity, $P$ is the pressure of the fluid,
$\phi$ is a force potential (e.g., centrifugal and gravitational forces),
$\chi(c)$ is the chemotactic sensitivity and $f(c)$ is the oxygen consumption
rate. The data $n_{0}\geq0,$ $c_{0}\geq0$ and $u_{0}$ with $\nabla\cdot
u_{0}=0$ are the initial conditions. A hypothesis for the model is
\begin{equation}
\chi(c),f(c),\chi^{\prime}(c),f^{\prime}(c)\geq0\text{ with }f(0)=0.
\label{Basic-cond-1}%
\end{equation}

System (\ref{CNSC}) was first introduced by Tuval \textit{et al.} \cite{Tuval}
and describes the motion of oxygen-driven swimming bacteria suspended in a
fluid (e.g., water) under the following basic effects: a chemotactic movement
towards oxygen (consumed by themselves); a convective transport of
microorganisms and substances mixed in the fluid; and a gravitational
influence due to the presence of the mass of a bacterial culture (heavier than
the fluid). For further details, see, for example, \cite{Lorz1}. Similar
models have also been analyzed in \cite{Pedley} and \cite{Hillesdon}.

In the case of constant decay rate $\gamma_{2}\geq0$ of a chemical $c$, and
without fluid coupling, a classical model related to (\ref{CNSC}) and
describing aggregation of microorganisms via a chemotaxis effect is the
so-called Keller--Segel system (see, e.g., \cite{Keller, Patlak})%

\begin{equation}%
\begin{cases}
\displaystyle\partial_{t}n=\Delta n-\nabla\cdot(n\chi\nabla c), & \qquad
x\in\mathbb{R}^{2}\ ,\ t>0,\ \\[3mm]%
\varepsilon\partial_{t}c-\Delta c=\gamma_{1}n-\gamma_{2}c, & \qquad
x\in\mathbb{R}^{2}\ ,\ t>0,\ \\[3mm]%
n(x,0)=n_{0}(x),\ \ c(x,0)=c_{0}(x),\ \  & \qquad x\in\mathbb{R}^{2},
\end{cases}
\label{KS}%
\end{equation}
where the chemical $c$ is produced by the microorganisms with the constant
rate $\gamma_{1}\geq0.$ The parameters $\varepsilon=0$ and $\varepsilon=1$
correspond respectively to the parabolic-elliptic and parabolic-parabolic
cases and the constant $\chi>0$ is the chemotactic sensitivity. System
\ref{KS} has been widely studied by means of different approaches and
functional settings, see, for example, the works \cite{Bedro-Masmoudi}%
,\cite{Biler1},\cite{Biler2},\cite{Blanchet},\cite{Calvez},\cite{Ferreira2}%
,\cite{Horstman0},\cite{Horstman1},\cite{Kozono1},\cite{Winkler1} and their
references. For $\varepsilon=1$ the global smoothness is an open problem and
global-in-time existence results are available under smallness conditions on
$n_{0},c_{0}$ in the corresponding functional setting. On the other hand, the
case $\varepsilon=0$ is completely resolved in the sense that there is a
threshold value $8\pi/\chi\gamma_{1}$ for the initial mass of $n_{0}$ that
decides between blow up ($\left\Vert n_{0}\right\Vert _{L^{1}}\geq8\pi
/\chi\gamma_{1}$) or not of solutions at a finite time $T,$ where the blow-up
nature is of mass-concentration such as Dirac delta measures (see, e.g.,
\cite{Blanchet-Carrillo1},\cite{Velazquez1}). For the case $\varepsilon=1,$
results about that kind of blow-up can be found in \cite{Nagai1, Nagai2}. The
system (\ref{KS}) in a bounded domain $\Omega\subset\mathbb{R}^{N}$ ($N\geq2$)
with the second equation replaced by $\partial_{t}c-\Delta c=-nc$ was analyzed
in \cite{Tao1}. There, considering homogeneous Neumann boundary conditions for
both $n$ and $c,$ the author developed a nice $L^{p}$-theory of global-in-time
solutions with data $u_{0}\in W^{1,r}(\Omega)\cap L^{1}(\Omega)$ and $c_{0}\in
W^{1,r}(\Omega),$ $r>n,$ under a smallness constraint on $\left\Vert
c_{0}\right\Vert _{L^{\infty}}$ depending on the dimension $N$ and the
chemotactic sensitivity $\chi$.

In view of possible concentration states, it is natural to wonder how the
system evolves from measures as initial densities, namely belonging to the
space of finite Radon measures $\mathcal{M}(\mathbb{R}^{2})$. In this
direction, Biler \cite{Biler1} showed local-in-time well-posedness for the
parabolic-parabolic variant of (\ref{KS}) with $n_{0}\in\mathcal{M}%
(\mathbb{R}^{2})$ and $c_{0}\in\dot{H}^{1}(\mathbb{R}^{2})$ such that
\[
\lim\sup_{t\rightarrow0^{+}}t^{1/3}\left\Vert e^{t\Delta}n_{0}\right\Vert
_{3/2}\,\ \text{is small enough,}%
\]
which imposes implicitly a smallness condition on the mass $\left\Vert
\cdot\right\Vert _{\mathcal{M}}$ of the singular part of $n_{0}$ (see
\cite[p.951]{Kato} and \cite{Giga1}). Moreover, global solutions are obtained
for $(n_{0},c_{0})$ sufficiently small in $\mathcal{M}(\mathbb{R}^{2}%
)\times\dot{H}^{1}(\mathbb{R}^{2}).$ With a slight adaptation, a version of
these results for $\varepsilon=0$ is also provided by \cite{Biler1}. After, in
line with the aforementioned threshold value for global existence, Bedrossian
and Masmoudi \cite{Bedro-Masmoudi} analyzed (\ref{KS}) with $\varepsilon
=\gamma_{2}=0$ and $\gamma_{1}=\chi=1$ and improved the smallness condition in
\cite{Biler1} by showing a local-in-time well-posedness result under the
condition $\max_{x\in\mathbb{R}^{2}}n_{0}(x)<8\pi$ on the atomic part of
$n_{0}.$

In comparison with (\ref{KS}), there are less results for the
chemotaxis-Navier-Stokes model (\ref{CNSC}). In the whole space $\mathbb{R}%
^{3},$ Duan \textit{et al.} \cite{Duan1} showed the global existence of
classical solutions for (\ref{CNSC}) with initial data being a small smooth
perturbation of the constant steady state $(n_{\infty},0,0)$ \ where
$n_{\infty}\geq0$ and $(n_{0},c_{0},u_{0})\rightarrow$ $(n_{\infty},0,0)$ as
$\left\vert x\right\vert \rightarrow\infty,$ as well \ as the convergence
rates of solutions towards $(n_{\infty},0,0).$ In the whole plane
$\mathbb{R}^{2}$, they analyzed a variant of (\ref{CNSC}) without the
nonlinearity $u\cdot\nabla u$ in the third equation (i.e., Stokes equations)
by considering smooth functions $\chi,f,\phi$ such that $\phi\geq0$,
$\nabla\phi\in L^{\infty}$ (and some finite weighted $L^{\infty}$-norms)$,$
$f^{\prime}>0$ and $f(0)=0$. Taking $n_{0}(1+\left\vert x\right\vert
+\left\vert \ln n_{0}\right\vert +\phi)\in L^{1}$, $u_{0}\in L^{2}$ and
$c_{0}\in L^{1}\cap L^{\infty}\cap H^{1}$ with small $\left\Vert
c_{0}\right\Vert _{L^{\infty}},$ they obtained existence of global weak
solutions. Another conditions on $c_{0},\chi,f$ considered by them were
$c_{0}\in L^{1}\cap L^{\infty}$, $\chi>0,$ $\chi^{\prime}\geq0,$ $f^{\prime
}>0,$ $f(0)=0$ and $(f/\chi)^{\prime\prime}<0$ with small $\left\Vert
c_{0}\right\Vert _{L^{4}}$ (and some small weighted $L^{\infty}$-norms of
$\nabla\phi$) instead of small $\left\Vert c_{0}\right\Vert _{L^{\infty}}.$

Chae \textit{et al.} \cite{Chae} considered the complete system (\ref{CNSC})
in $\mathbb{R}^{2}$ and proved the existence-uniqueness of classical solutions
by assuming
\begin{equation}
(n_{0},c_{0},u_{0})\in H^{m-1}\times H^{m}\times H^{m}\text{ with }%
m\geq3\text{ and }\nabla^{l}\phi\in L^{\infty}\text{ for all }1\leq\left\vert
l\right\vert \leq m, \label{Cond-intro-2}%
\end{equation}
$\phi\geq0,$ and supposing that $\chi,f\in C^{m}(\mathbb{R}^{+})$,
(\ref{Basic-cond-1}) and that there exists a constant $\mu>0$ such that
\begin{equation}
\sup_{c}|\chi(c)-\mu f(c)|<\epsilon\text{ for sufficiently small }\epsilon>0.
\label{Cond-intro-2-1}%
\end{equation}
Still in the same work, by considering $\chi,f\in C^{1}(\mathbb{R}^{+})$
satisfying (\ref{Basic-cond-1}) and $\chi(c)=\mu f(c)$ for some constant
$\mu>0$, and taking $\phi\geq0$ with $\nabla^{l}\phi\in L^{\infty}$ for all
$1\leq\left\vert l\right\vert \leq2,$ they treated the 3D case and obtained
the existence of a global weak solution with $n_{0}\in L^{1},$ $c_{0}\in
H^{1}\cap L^{\infty},$ $u_{0}\in L^{2}$ and $n_{0}(1+\left\vert x\right\vert
+\left\vert \ln n_{0}\right\vert ))\in L^{1}.$ In \cite{Chae2} the same
authors of \cite{Chae} proved the existence of a unique classical solutions
for (\ref{CNSC}) in $\mathbb{R}^{2}$ by considering (\ref{Cond-intro-2}),
$\chi,f\in C^{m}(\mathbb{R}^{+})$, and (\ref{Basic-cond-1}), and supposing a
smallness condition only on $\left\Vert c_{0}\right\Vert _{L^{\infty}}$ (see
also \cite{Chae3}). Assuming (\ref{Cond-intro-2}), $n_{0}(1+\left\vert
x\right\vert +\left\vert \ln n_{0}\right\vert ))\in L^{1}$ and $c_{0}\in
L^{1},$ and considering $\chi,f\in C^{m}(\mathbb{R})$ with $f\geq0$ and
$f(0)=0$, Duan \textit{et al.} \cite{Duan2} obtained a unique global classical
solution in $\mathbb{R}^{2}$ under a smallness condition involving $\left\Vert
c_{0}\right\Vert _{L^{\infty}},$ $\left\Vert n_{0}\right\Vert _{L^{1}}%
,\chi,f,$ and a Gagliardo-Nirenberg type constant, which is particularly
verified for $\left\Vert c_{0}\right\Vert _{L^{\infty}}$ or $\left\Vert
n_{0}\right\Vert _{L^{1}}$ small enough. They also analyzed the case of a
bounded domain $\Omega\subset$ $\mathbb{R}^{2}$ by assuming a similar
condition on $\left\Vert c_{0}\right\Vert _{L^{\infty}},$ $\left\Vert
n_{0}\right\Vert _{L^{1}},\chi,f$ but that imposes necessarily a restriction
on the size of $\left\Vert n_{0}\right\Vert _{L^{1}}.$

In the case of bounded convex domains $\Omega\subset\mathbb{R}^{2}$, the
existence-uniqueness of global classical solutions was proved by Winkler
\cite{Winkler2} via suitable entropy estimates and without assuming smallness
conditions on the initial data. He considered $\phi\in C^{2}(\overline{\Omega
})$, a class of bounded continuous initial data $(n_{0},c_{0},u_{0})$ with
$n_{0},c_{0}>0$ in $\overline{\Omega}$, and $\chi,f\in C^{2}([0,\infty))$ with
$\chi>0$ in $[0,\infty)$ $,$ $f>0$ in $(0,\infty),$ $f(0)=0$ such that
\begin{equation}
\left(  \frac{f}{\chi}\right)  ^{\prime}>0,\text{ }\left(  \frac{f}{\chi
}\right)  ^{\prime\prime}\leq0\text{ and }(\chi\cdot f)^{\prime}\geq0.
\label{Cond-intro-3}%
\end{equation}
The conditions in \cite{Winkler2} (see also \cite{Chae2,Duan2}) allow to cover
a relevant case, namely the so-called prototypical choice
\begin{equation}
\chi(c)\equiv\kappa_{1}\text{ and }f(c)=\kappa_{2}c, \label{prototypical1}%
\end{equation}
that is, the chemotactic sensitivity is constant and the oxygen consumption
rate is proportional to the oxygen concentration itself (see \cite{Rosen1}).
In fact, this can be seen as a borderline case in \cite{Winkler2}. For
simplicity, we assume $\kappa_{1}=\kappa_{2}=1$ in (\ref{prototypical1}). For
further results in bounded domains of $\mathbb{R}^{2}$ and $\mathbb{R}^{3},$
the reader is refereed to \cite{Liu-Lorz},\cite{Lorz1},\cite{Winkler2}%
,\cite{Winkler3} and references therein.

Motivated by \cite{Winkler2}, Zhang and Zheng \cite{Zhang-Zheng} approached
the prototypical case (\ref{prototypical1}) for (\ref{CNSC}) in the whole
$\mathbb{R}^{2}$ by developing a unified approach for weak solutions and more
general initial data. More precisely, by means of a scale decomposition
technique and Zygmund spaces, and considering $\nabla\phi\in L^{\infty}$ and
the class of rough initial data
\begin{equation}
\left\{  n_{0}\in L^{1}\cap L^{2},\text{ }\nabla\sqrt{c_{0}}\in L^{2},\text{
}c_{0}\in L^{1}\cap L^{\infty},\text{ }u_{0}\in L^{2}\right\}
,\label{class-Z}%
\end{equation}
they obtained the existence-uniqueness of global distributional solutions in
the class
\begin{align*}
n &  \in L^{\infty}(\mathbb{R}^{+};L^{1})\cap L_{loc}^{\infty}(\mathbb{R}%
^{+};L^{2})\cap L_{loc}^{2}(\mathbb{R}^{+};H^{1}),\\
c &  \in L^{\infty}(\mathbb{R}^{+};L^{1}\cap L^{\infty})\cap L_{loc}^{\infty
}(\mathbb{R}^{+};H^{1})\cap L_{loc}^{2}(\mathbb{R}^{+};H^{2}),\\
u &  \in L_{loc}^{\infty}(\mathbb{R}^{+};L^{2})\cap L_{loc}^{2}(\mathbb{R}%
^{+};H^{1}).
\end{align*}
As a matter of fact, in the whole $\mathbb{R}^{2}$ the problem of formation of
singularity in finite time for (\ref{CNSC})-(\ref{Basic-cond-1}) seems to be
still open (see \cite[pp. 352-353]{Chae3}), even in the prototypical case. In
view of the behavior of system (\ref{KS}) in the parabolic-elliptic case with
$\gamma_{2}=0$, a possible singularity should probably present a scenario with
mass and/or vortex concentrations corresponding to a solution taking a
singular measure value at the blow-up time. Such concentrations could occur in
points or curves, which naturally motivates studying how (\ref{CNSC}) evolves
from these singular states and whether solutions could be continued in
suitable measure spaces. In any case, from a physical/biological point of
view, understanding the model under very intense concentrations could provide
a test for its robustness, particularly for its behavior in a scenario of
extremely high concentration of bacteria and for the presence of vortex-sheets.

In this paper our aim is to provide an existence-uniqueness theory for system
(\ref{CNSC}) evolving from arbitrary initial Radon measures in the pivotal
prototypical case (\ref{prototypical1}). More precisely, we obtain
global-in-time solutions with initial density $n_{0}\in\mathcal{M}%
(\mathbb{R}^{2})$ and the fluid vorticity $\zeta_{0}\in\mathcal{M}%
(\mathbb{R}^{2})$ regardless of their sizes, see Theorems \ref{TLS} and
\ref{TLSM} for the local-in-time theory, and Theorems \ref{global} and
\ref{global-M} for the global one. In principle, our solutions are of mild
type, but they smooth out instantaneously for $t>0$ (see Remark
\ref{regularity} and Remarks \ref{rem-teo-measure} $(iii)$ and
\ref{Rem-global} $(ii)$) and thereby verify classically the system for
$x\in\mathbb{R}^{2}$ and $t>0$. For the oxygen concentration $c_{0}$, we
assume $c_{0}\in L^{\infty}(\mathbb{R}^{2})$ with $\nabla c_{0}\in
L^{2}(\mathbb{R}^{2})$ and a smallness condition on $\left\Vert c_{0}%
\right\Vert _{L^{\infty}}.$ This condition can be removed in the local theory
(see Remark \ref{RTL} $(i)$ and $(ii)$). The force $\nabla\phi$ is taken in
$H^{2}(\mathbb{R}^{2})$ for global solutions and (even locally) in the case of
singular measures $n_{0},\zeta_{0}$. However, for local-in-time solutions, we
can assume only $\nabla\phi\in L^{2}(\mathbb{R}^{2})$ when $n_{0},\zeta_{0}$
are absolutely continuous measures. In comparison with the previous works, our
results provide a new initial-data class for $n_{0}$ and $\zeta_{0}$, covering
particularly singular measures such as Dirac delta, filaments in
$\mathbb{R}^{2}$ (measures concentrated on smooth curves; vortex-sheets for
$\zeta$), among others. Also, our conditions on $\phi$ are different from
previous results on mild or classical solutions for (\ref{CNSC}) in the whole
$\mathbb{R}^{2}$, see (\ref{Cond-intro-2}) and references \cite{Chae}%
,\cite{Chae2},\cite{Chae3},\cite{Duan1},\cite{Duan2}. The local-in-time
uniqueness, which yields the global one (by standard arguments), is discussed
in details in item $(ii)$ of Remark \ref{rem-teo-measure}. So, we have the
uniqueness of solution when $n_{0},\zeta_{0}\in L^{1}(\mathbb{R}^{2})$ or,
more generally, $n_{0},\zeta_{0}\in\mathcal{M}(\mathbb{R}^{2})$ have small
atomic part (e.g., non-atomic measures), see Theorem \ref{global-M}. It is
worthy to recall we are always considering nonnegative data $n_{0},c_{0}\geq0$
and then the obtained solutions $(n,c,\zeta)$ satisfy $n,c\geq0$, because the
equations for $n,c$ preserve the initial-data sign.

Inspired by previous works for the 2D Navier-Stokes equations treating with
initial vorticity measures, such as \cite{Bedro-Masmoudi},\cite{Kato}%
,\cite{Giga1}, we consider (\ref{CNSC})-(\ref{prototypical1}) with the
vorticity formulation for the fluid movement. Then, computing the vorticity of
the fluid $\zeta=\nabla^{\perp}\cdot u=\partial_{1}u_{2}-\partial_{2}u_{1}$ in
the third equation of (\ref{CNSC}) and considering (\ref{prototypical1}), we
are led to its vorticity formulation%

\begin{equation}%
\begin{cases}
\displaystyle\partial_{t}n+u\cdot\nabla n=\Delta n-\nabla\cdot(n\nabla c), &
\qquad x\in\mathbb{R}^{2}\ ,\ t>0,\ \\[3mm]%
\partial_{t}c+u\cdot\nabla c=\Delta c-cn, & \qquad x\in\mathbb{R}%
^{2}\ ,\ t>0,\ \\[3mm]%
\partial_{t}\zeta+\nabla\cdot(\zeta S\ast\zeta)=\Delta\zeta+\nabla^{\perp
}\cdot(n\nabla\phi), & \qquad x\in\mathbb{R}^{2}\ ,\ t>0,\ \\[3mm]%
n(x,0)=n_{0},\ \ c(x,0)=c_{0},\ \ \zeta(x,0)=\zeta_{0}, & \qquad
x\in\mathbb{R}^{2},
\end{cases}
\label{CNS}%
\end{equation}
where the velocity can be recovered via the Biot-Savart law
\begin{equation}
u=S\ast\zeta\text{ with }S(x)=(2\pi)^{-1}|x|^{-2}(x_{2},-x_{1})\text{ and
}x=(x_{1},x_{2}). \label{velocity1}%
\end{equation}
In (\ref{velocity1}), the symbol $\ast$ stands for the convolution operator
and $S(x)$ is the so-called Biot-Savart kernel.

In our analysis of (\ref{CNS}), we employ some Kato-type classes that are
time-weighted ones based on $L^{p}$-spaces (see (\ref{spacesolutions})), as
well as obtain suitable $L^{p}$-estimates uniformly in time (see, e.g.,
Propositions \ref{zetaboundedness} and \ref{GCboundedness}). Here we are also
inspired by some arguments found in \cite{Tao1}. Moreover, we consider
critical cases of those classes, that is, a setting invariant under the
scaling map $(n,c,\zeta)\rightarrow(n_{\lambda},c_{\lambda},\zeta_{\lambda})$
for each $\lambda>0$, where the rescaled triple $(n_{\lambda},c_{\lambda
},\zeta_{\lambda})$ is given by%

\begin{equation}
(n_{\lambda}(x,t),c_{\lambda}(x,t),\zeta_{\lambda}(x,t)):=(\lambda
^{2}n(\lambda x,\lambda^{2}t),c(\lambda x,\lambda^{2}t),\lambda^{2}%
\zeta(\lambda x,\lambda^{2}t)). \label{scalling}%
\end{equation}

This paper is organized as follows. In Section \ref{sectionP} we give some
preliminaries about Kato-type spaces, heat semigroup and Biot-Savart operator,
and recall a useful fixed-point lemma. In Section \ref{sectionL}, we obtain
the local-in-time existence and uniqueness of solutions for (\ref{CNS}) with
initial data $n_{0},\zeta_{0}\in L^{1}(\mathbb{R}^{2})$ and $c_{0}\in
L^{\infty}(\mathbb{R}^{2}),$ as well as give some further properties. Section
\ref{sectionM} is devoted to prove the local existence of solutions for
(\ref{CNS}) with Radon measures $\mathcal{M}(\mathbb{R}^{2})$ by extending the
$L^{1}$-result of Section \ref{sectionL}. Finally, in Section \ref{sectionG}
we show that the previously obtained solutions can be extended globally in time.


\section{Preliminaries}

\label{sectionP}

This section is devoted to give some notations and recall basic properties
about the heat semigroup and the Biot-Savart operator, as well as an abstract
fixed-point lemma.

We begin by defining suitable time-weighted spaces, namely the so-called
Kato-type spaces (see, e.g., \cite{Kato}), which will be used in the analysis
of (\ref{CNS}). Given a Banach space $X$ with norm $||\cdot||_{X}$ and an
interval $I\subset\mathbb{R}$, we denote by $C(I;X)$ the set of all continuous
$f:I\rightarrow X.$ The space $BC(I;X)$ consists of all bounded functions in
$C(I;X)$. For $\alpha\geq0$ and $T\in(0,\infty]$, the Kato space $C_{\alpha
}((0,T);X)$ stands for the set of all continuous functions $f:(0,T)\rightarrow
X$ such that%
\[
{|||}f{|||}_{X,\alpha}:=\sup_{0<t<T}t^{\alpha}||f(t)||_{X}.
\]
The space $C_{\alpha}((0,T);X)$ is a Banach space with the norm ${|||}%
\cdot{|||}_{X,\alpha}$.

\begin{remark}
\label{space-norm1}For simplicity, we often omit the interval $(0,T)$, i.e.
write only $C_{\alpha}(X)$. Also, for $X=L^{p}(\mathbb{R}^{2})$, we write
$C_{\alpha}(L^{p})$ and ${|||}\cdot{|||}_{p,\alpha}$ instead of $C_{\alpha
}(L^{p}(\mathbb{R}^{2}))$ and ${|||}\cdot{|||}_{L^{p}(\mathbb{R}^{2})}$,
respectively. Moreover, we sometimes consider the natural adaptation of those
definitions for intervals $(T_{0},T_{1})$, with $T_{0}>0$. In this case, the
norm in $C_{\alpha}((T_{0},T_{1});X)$ is given by
\end{remark}

\[
{|||}f{|||}_{X,\alpha}=\sup_{T_{0}<t<T_{1}}(t-T_{0})^{\alpha}||f(t)||_{X}.
\]

In $C^{\alpha}((0,T);X)$, we can define the seminorm%

\begin{equation}
\pmb{\lVert}f\pmb{\lVert}_{X,\alpha}:=\limsup_{t\rightarrow0}t^{\alpha
}||f(t)||_{X}. \label{seminorm1}%
\end{equation}

By definition, the value of $\pmb{\lVert}f\pmb{\lVert}_{X,\alpha}$ depends
only on the values of $f(t)$ when $t$ is close to $0$. Note that
$\pmb{\lVert}f\pmb{\lVert}_{X,\alpha}\leq{|||}f{|||}_{X,\alpha},$ for all
$f\in C^{\alpha}((0,T);X)$. For $X=L^{p},$ we denote
$\pmb{\lVert}f\pmb{\lVert}_{p,\alpha}=\pmb{\lVert}f\pmb{\lVert}_{L^{p},\alpha
}$. The subset of $C_{\alpha}((0,T);X)$ consisting of all $f$ such that
$\pmb{\lVert}f\pmb{\lVert}_{X,\alpha}=0$ will be denoted by $\dot{C}_{\alpha
}((0,T);X)$.

Let $\{e^{t\Delta}\}_{t\geq0}$ denote the heat semigroup. For any
$f\in\mathcal{S^{\prime}}$, it is well known that $e^{t\Delta}$ is the
convolution operator
\begin{equation}
e^{t\Delta}f=h(\cdot,t)\ast f, \label{heat-1}%
\end{equation}
where $h(x,t)$ is the so-called heat kernel $h(x,t)=(4\pi t)^{-1}\exp
(-|x|^{2}/4t)$. In the $L^{p}$-setting, using Young inequality, we have the
estimate
\begin{equation}
\left\Vert \nabla_{x}^{m}e^{t\Delta}f\right\Vert _{q}\leq c_{\gamma}%
t^{-\frac{\left\vert m\right\vert }{2}-(\frac{1}{r}-\frac{1}{q})}||f||_{r},
\label{est-sg-1}%
\end{equation}
for all $f\in L^{r}(\mathbb{R}^{2})$, where $1\leq r\leq q<\infty,$
$m\in(\mathbb{N}\cup\{0\})^{2},$ and $c_{\gamma}>0$ is a constant depending
only on $\gamma=\frac{1}{r}-\frac{1}{q}$ and $m$. In view of (\ref{est-sg-1}),
we have that $e^{t\Delta}f\in C_{\alpha}(L^{p}(\mathbb{R}^{2}))$ provided that
$f\in L^{r}(\mathbb{R}^{2})$ with the estimate%
\[
{|||}e^{t\Delta}f{|||}_{q,\gamma}\leq c_{\gamma}||f||_{r}.
\]
Moreover, if $r<q$, then $\gamma>0$ and
\begin{equation}
e^{t\Delta}f\in\dot{C}_{\gamma}(L^{q}). \label{Prop-point-1}%
\end{equation}

Let $\mathcal{M}(\mathbb{R}^{n})$ be the space of Radon measures in
$\mathbb{R}^{n}$ endowed with the total variation norm $\pmb{|}\mu\pmb{|}.$ We
recall the Young inequality
\begin{equation}
||f\ast\mu||_{p}\leq||f||_{p}\pmb{|}\mu\pmb{|}, \label{young-2}%
\end{equation}
where $f\in L^{p}(\mathbb{R}^{n})$ $(1\leq p\leq\infty)$ and $\mu
\in\mathcal{M}(\mathbb{R}^{n})$. For $n=2,$ applying (\ref{young-2}) in
(\ref{heat-1}), we arrive at the heat estimate
\begin{equation}
{|||}e^{t\Delta}\omega{|||}_{q,1-1/q}\leq c_{q}\pmb{|}\omega\pmb{|},
\label{1.7}%
\end{equation}
for all $\omega\in\mathcal{M}(\mathbb{R}^{2})$, where $q\in\lbrack1,\infty]$
and the constant $c_{q}>0$ depends only on $q$.

In the existence proof of local-in-time solutions, we are going to employ a
contraction argument. In order to avoid exhaustive fixed-point computations,
the following abstract lemma will be useful. For the proof, see, for example,
\cite[Lemma 4.1]{Ferreira1}.

\begin{lemma}
\label{LA} For $1\leq i\leq3$, let $X_{i}$ be a Banach space with norm
$||\cdot||_{x_{i}}$. Consider the Banach space $X=X_{1}\times X_{2}\times
X_{3}$ endowed with the norm
\[
||x||_{X}=||x_{1}||_{X_{1}}+||x_{2}||_{X_{2}}+||x_{3}||_{X_{3}},
\]
where $x=(x_{1},x_{2},x_{3})\in X$. For $1\leq i,j,k\leq3$, assume that
$B_{i,j}^{k}:X_{i}\times X_{j}\rightarrow X_{k}$ is a continuous bilinear map,
i.e., there is a constant $C_{i,j}^{k}$ (depending only on $k,i,j$) such that
\[
||B_{i,j}^{k}(x_{i},x_{j})||_{X_{k}}\leq C_{i,j}^{k}||x_{i}||_{X_{i}}%
||x_{j}||_{X_{j}},\hspace{1cm}\forall x_{i}\in X_{i},\ \ x_{j}\in X_{j}.
\]

Assume also that $L_{1}^{3}:X_{1}\rightarrow X_{3}$ is a continuous linear
map, i.e., there is a constant $\alpha$ such that%
\[
||L_{1}^{3}(x_{1})||_{X_{3}}\leq\alpha||x_{1}||_{X_{1}},\hspace{1cm}\forall
x_{1}\in X_{1}.
\]

\end{lemma}

Set $K_{1}:=1+\alpha$ and $K_{2}:=\alpha\sum_{i,j=1}^{3}C_{i,j}^{1}%
+\sum_{k,i,j=1}^{3}C_{i,j}^{k}$. Take $0<\epsilon<\frac{1}{4K_{1}K_{2}}$ and
define $\mathcal{B}_{\epsilon}:=\{x\in X:||x||_{X}\leq2K_{1}\epsilon\}$. If
$||y||_{X}\leq\epsilon$, then there exists a unique solution $x\in
\mathcal{B}_{\epsilon}$ for the equation $x=y+B(x),$ where $y=(y_{1}%
,y_{2},y_{3})$ and $B(x)=(B_{1}(x),B_{2}(x),B_{3}(x))$ with
\[
B_{1}(x)=\sum_{i,j=1}^{3}B_{i,j}^{1}(x_{i},x_{j}),\text{ \ }B_{2}%
(x)=\sum_{i,j=1}^{3}B_{i,j}^{2}(x_{i},x_{j}),
\]
and%
\[
B_{3}(x)=\sum_{i,j=1}^{3}B_{i,j}^{3}(x_{i},x_{j})+(L_{1}^{3}\circ(y_{1}%
+B_{1}))(x).
\]

Finally, in view of the Hardy-Littlewood-Sobolev theorem for fractional
integration (see \cite[p.3]{Grafakos}), we have that the linear operator
$S\ast$ defined in (\ref{velocity1}) can be handled in the $L^{p}$-setting as
follows
\begin{equation}
||S\ast\psi||_{p}\leq\sigma_{q}||\psi||_{q},\ \label{littlewood}%
\end{equation}
where $q\in(1,2),$ $\frac{1}{p}=\frac{1}{q}-\frac{1}{2},$ and $\sigma_{q}$ is
a constant depending only on $q$.

\section{Local-in-time solutions with $L^{p}$-data}

\label{sectionL} Consider the space
\begin{equation}
X(a,T)=X_{1}(a,T)\times X_{2}(a,T)\times X_{3}(a,T), \label{spacesolutions-0}%
\end{equation}
where $0\leq a<T\leq\infty$ and
\begin{equation}%
\begin{split}
X_{1}(a,T)=  &  C_{\alpha_{1}}((a,T);L^{p_{1}}(\mathbb{R}^{2})),\\
X_{2}(a,T)=  &  \{c\in C_{0}((a,T);L^{\infty}(\mathbb{R}^{2}):\nabla
c\in{C_{\alpha_{2}}((a,T);}L^{p_{2}}(\mathbb{R}^{2})\},\\
X_{3}(a,T)=  &  C_{\alpha_{3}}((a,T);L^{p_{3}}(\mathbb{R}^{2})).
\end{split}
\label{spacesolutions}%
\end{equation}
For simplicity, we will often omit the interval $(a,T)$ in $X(a,T)$ and
$X_{i}(a,T)$ by writing only $X$ and $X_{i}$, respectively. The norms in
$X_{1},$ $X_{2}$ and $X_{3}$ are given respectively by
\begin{equation}%
\begin{split}
||n||_{X_{1}}:=  &  |||n|||_{p_{1},\alpha_{1}}=\sup_{a<t<T}(t-a)^{\alpha_{1}%
}||n(t)||_{p_{1}},\\
||c||_{X_{2}}:=  &  |||c|||_{\infty,0}+|||\nabla c|||_{p_{2},\alpha_{2}}%
=\sup_{a<t<T}||c(t)||_{\infty}+\sup_{a<t<T}(t-a)^{\alpha_{2}}||\nabla
c(t)||_{p_{2}},\\
||\zeta||_{X_{3}}:=  &  |||\zeta|||_{p_{3},\alpha_{3}}=\sup_{a<t<T}%
(t-a)^{\alpha_{3}}||\zeta(t)||_{p_{3}},
\end{split}
\label{norm}%
\end{equation}
and $||g||_{X}=||n||_{X_{1}}+||c||_{X_{2}}+||\zeta||_{X_{3}}$, for all
$g=(n,c,\zeta)\in X$.

According to (\ref{scalling}), in order to consider critical spaces, we assume
that $\alpha_{i}$ and $p_{i}$ satisfy
\begin{equation}
\alpha_{1}+\frac{1}{p_{1}}=1,\ \ \alpha_{2}+\frac{1}{p_{2}}=\frac{1}%
{2},\ \ \alpha_{3}+\frac{1}{p_{3}}=1. \label{critical}%
\end{equation}
Back to system (\ref{CNS}), one can formally convert it into the integral
formulation (via Duhamel principle)%

\begin{equation}%
\begin{cases}
\displaystyle n(t)=e^{t\Delta}n_{0}-\int_{0}^{t}\nabla\cdot e^{(t-s)\Delta
}(n\nabla c)(s)ds-\int_{0}^{t}\nabla\cdot e^{(t-s)\Delta}(n(S\ast
\zeta))(s)ds,\\[3mm]%
\displaystyle c(t)=e^{t\Delta}c_{0}-\int_{0}^{t}e^{(t-s)\Delta}(nc)(s)ds-\int
_{0}^{t}e^{(t-s)\Delta}((S\ast\zeta)\cdot\nabla c)(s)ds,\\[3mm]%
\displaystyle\zeta(t)=e^{t\Delta}\zeta_{0}-\int_{0}^{t}\nabla\cdot
e^{(t-s)\Delta}(\zeta S\ast\zeta)(s)ds-\int_{0}^{t}\nabla^{\perp}\cdot
e^{(t-s)\Delta}(n\nabla\phi)(s)ds.
\end{cases}
\label{IF}%
\end{equation}

A triple $(n,c,\zeta)$ satisfying (\ref{IF}) is called a \textit{mild}
solution of (\ref{CNS}). In order to simplify the notation, we define the
following operators linked to (\ref{IF})
\begin{equation}
B_{1,2}^{1}(n,c)(t):=\int_{0}^{t}\nabla\cdot e^{(t-s)\Delta}(n\nabla
c)(s)ds,\ \ \ \ \ \ \ B_{1,3}^{1}(n,\zeta)(t):=\int_{0}^{t}\nabla\cdot
e^{(t-s)\Delta}(n(S\ast\zeta))(s)ds, \label{Op1}%
\end{equation}

\begin{equation}
B_{2,3}^{2}(c,\zeta)(t):=\int_{0}^{t}e^{(t-s)\Delta}((S\ast\zeta)\cdot\nabla
c)(s)ds,\ \ \ \ \ \ \ B_{1,2}^{2}(n,c)(t):=\int_{0}^{t}e^{(t-s)\Delta
}(nc)(s)ds, \label{Op2}%
\end{equation}

\begin{equation}
B_{3,3}^{3}(\zeta,\tilde{\zeta})(t):=\int_{0}^{t}\nabla\cdot e^{(t-s)\Delta
}(\zeta S\ast\tilde{\zeta})(s)ds,\ \ \ \ \ \ \ L_{1}^{3}(n)(t):=\int_{0}%
^{t}\nabla^{\perp}\cdot e^{(t-s)\Delta}(n\nabla\phi)(s)ds. \label{Op3}%
\end{equation}
Note that $B_{i,j}^{k}$ are bilinear operators while $L_{1}^{3}$ is linear.

Adding this notation to the integral formulation (\ref{IF}), we obtain the
more compact expression%

\begin{equation}%
\begin{cases}
\displaystyle n(t)=e^{t\Delta}n_{0}-B_{1,2}^{1}(n,c)(t)-B_{1,3}^{1}%
(n,\zeta)(t),\\[3mm]%
\displaystyle c(t)=e^{t\Delta}c_{0}-B_{2,3}^{2}(c,\zeta)(t)-B_{1,2}%
^{2}(n,c)(t),\\[3mm]%
\displaystyle\zeta(t)=e^{t\Delta}\zeta_{0}-B_{3,3}^{3}(\zeta,{\zeta}%
)(t)-L_{1}^{3}(n)(t).
\end{cases}
\label{IFN}%
\end{equation}

Next we set%

\begin{equation}%
\begin{split}
\Phi(n,c,\zeta)(t)  &  =(e^{t\Delta}n_{0}-B_{1,2}^{1}(n,c)(t)-B_{1,3}%
^{1}(n,\zeta)(t),e^{t\Delta}c_{0}-B_{2,3}^{2}(c,\zeta)(t)-B_{1,2}%
^{2}(n,c)(t),\\
&  e^{t\Delta}\zeta_{0}-B_{3,3}^{3}(\zeta,{\zeta})(t)-L_{1}^{3}(n)(t)).
\end{split}
\label{Op-fi}%
\end{equation}
Note that a mild solution of system (\ref{CNS}) can be seen as a fixed point
of $\Phi$. Back to the functional setting for solutions of (\ref{CNS}), we
need that the index $p_{i}$ and $\alpha_{i}$ in (\ref{spacesolutions}) satisfy
some technical conditions in order to obtain suitable estimates for the
operators in (\ref{IFN}). For the sake of simplicity, we will call it
\textbf{Condition A}.

\begin{definition}
\label{Cond-A}For $1\leq i\leq3$, let $p_{i}\in\lbrack1,\infty]$ and
$\alpha_{i}\geq0$. We say that $p_{i},\alpha_{i}$ satisfies the
\textbf{Condition A} if $p_{i}$ and $\alpha_{i}$ satisfies
\begin{equation}
\alpha_{1}+\frac{1}{p_{1}}=1;\ \ \alpha_{2}+\frac{1}{p_{2}}=\frac{1}%
{2};\ \ \alpha_{3}+\frac{1}{p_{3}}=1; \label{A1}%
\end{equation}%
\begin{equation}
\frac{1}{p_{1}}+\frac{1}{p_{2}}\leq1;\ \ p_{2}\in(2,\infty];\ \ p_{3}%
\in\lbrack4/3,2);\ \ \frac{1}{p_{1}}+\frac{1}{p_{3}}\leq\frac{3}{2};
\label{A2}%
\end{equation}%
\begin{equation}
\frac{1}{p_{2}}+\frac{1}{p_{3}}<\frac{3}{2};\ \ p_{1}\geq2;\ \ 0\leq\frac
{1}{p_{1}}-\frac{1}{p_{2}}<\frac{1}{2};\ \ 0<\frac{1}{p_{3}} - \frac{1}{p_{1}%
}\leq\frac{1}{2} ; \label{A3}%
\end{equation}%
\begin{equation}
\alpha_{1}+\alpha_{2}<1;\ \ \alpha_{1}+\alpha_{3}<1;\ \ \alpha_{2}+\alpha
_{3}<1;\ \ \alpha_{3}<\frac{1}{2}. \label{A4}%
\end{equation}

\end{definition}

\begin{remark}
\label{Rem-Cond-A}Note that the above condition is not empty. For example, the
index set $p_{1}=\frac{17}{8}$, $p_{2}=3$, $p_{3}=\frac{15}{8}$, $\alpha
_{1}=\frac{9}{17}$, $\alpha_{2}=\frac{1}{6}$ and $\alpha_{3}=\frac{7}{15}$
satisfies \textbf{Condition A}.
\end{remark}

Our first existence and uniqueness result for system (\ref{CNS}) reads as follows.

\begin{theorem}
{(Local-in-time solutions)} \label{TLS} For $1\leq i\leq3$, let $p_{i}%
,\alpha_{i}$ satisfy \textbf{Condition A} in Definition \ref{Cond-A}. Let
$\zeta_{0},n_{0}\in L^{1}(\mathbb{R}^{2})$, $c_{0}\in L^{\infty}%
(\mathbb{R}^{2})$ with $\nabla c_{0}\in L^{2}(\mathbb{R}^{2})$, and
$\nabla\phi\in L^{2}(\mathbb{R}^{2})$. If $||c_{0}||_{\infty}$ is sufficiently
small, then there exists a $T_{0}>0$ such that system (\ref{CNS}) has a unique
mild solution $(n,c,\zeta)\in X$, locally in $(0,T_{0})$, with initial data
$(n_{0},c_{0},\zeta_{0})$.
\end{theorem}

\textit{Proof.} We start by showing that the operators $B_{i,j}^{k}$ defined
in (\ref{Op1})-(\ref{Op3}) are continuous bilinear map, with the bilinear
constants $C_{i,j}^{k}$ independent of the size of the time interval $(0,T)$.

For $B_{1,2}^{1}$, using (\ref{est-sg-1}), $\frac{1}{p_{1}}+\frac{1}{p_{2}%
}\leq1$, and H\"{o}lder inequality, we can estimate%

\begin{align}
||B_{1,2}^{1}(n,c)(t)||_{p_{1}}  &  =\bigg|\bigg|\int_{0}^{t}\nabla\cdot
e^{(t-s)\Delta}(n\nabla c)(s)ds\bigg|\bigg|_{p_{1}}\nonumber\\
&  \leq\int_{0}^{t}||\nabla\cdot e^{(t-s)\Delta}(n\nabla c)(s)||_{p_{1}%
}ds\nonumber\\
&  \leq C\int_{0}^{t}{(t-s)^{-\frac{1}{2}-(\frac{1}{q}-\frac{1}{p_{1}})}%
}||(n\nabla c)(s)||_{q}ds,\ \ \ \ \text{with}\ \ \frac{1}{q}=\frac{1}{p_{1}%
}+\frac{1}{p_{2}}\nonumber\\
&  \leq C\int_{0}^{t}{(t-s)^{-\frac{1}{2}-\frac{1}{p_{2}}}}||n(s)||_{p_{1}%
}||\nabla c(s)||_{p_{2}}ds.\ \ \ \ \label{aux-B112-1}%
\end{align}
Recalling the notation of the norm $|||\cdot|||_{p_{1},\alpha_{1}}$ given in
Remark \ref{space-norm1} and Beta function
\begin{equation}
\beta(\delta_{1},\delta_{2})=\int_{0}^{1}t^{\delta_{1}-1}(1-t)^{\delta_{2}%
-1}dt, \label{beta}%
\end{equation}
and employing the change $s\rightarrow st$ and the fact that $\alpha_{2}%
+\frac{1}{p_{2}}=\frac{1}{2}$ (see (\ref{A1})), we obtain that
\begin{align}
\text{R.H.S. of (\ref{aux-B112-1})}  &  \leq C|||n|||_{p_{1},\alpha_{1}%
}|||\nabla c|||_{p_{2},\alpha_{2}}\int_{0}^{t}{(t-s)^{-\frac{1}{2}-\frac
{1}{p_{2}}}}s^{-\alpha_{1}-\alpha_{2}}ds\nonumber\\
&  =C|||n|||_{p_{2},\alpha_{1}}|||\nabla c|||_{p_{2},\alpha_{2}}t^{\frac{1}%
{2}-\frac{1}{p_{2}}-\alpha_{1}-\alpha_{2}}\beta\bigg(\frac{1}{2}-\frac
{1}{p_{2}},1-\alpha_{1}-\alpha_{2}\bigg)\ \ \ \ \nonumber\\
&  \leq C|||n|||_{p_{1},\alpha_{1}}|||\nabla c|||_{p_{2},\alpha_{2}}%
t^{-\alpha_{1}}\beta\bigg(\frac{1}{2}-\frac{1}{p_{2}},1-\alpha_{1}-\alpha
_{2}\bigg)\ \ \text{by}\ \text{\ (\ref{A1})}\nonumber\\
&  \leq C||n||_{X_{1}}||c||_{X_{2}}t^{-\alpha_{1}}\beta\bigg(\frac{1}{2}%
-\frac{1}{p_{2}},1-\alpha_{1}-\alpha_{2}\bigg)\text{,}\ \label{aux-B112-2}%
\end{align}
where above we also used the norms defined in (\ref{norm}) and the conditions
$p_{2}>2$ and (\ref{A4}). Inserting (\ref{aux-B112-2}) into (\ref{aux-B112-1})
and afterwards multiplying both sides of the resulting inequality by
$t^{\alpha_{1}}$ and computing the supremum over $t$, we arrive at
\[
||B_{1,2}^{1}(n,c)||_{X_{1}}\leq C_{1,2}^{1}||n||_{X_{1}}||c||_{X_{2}},
\]
where $C_{1,2}^{1}:=C\beta\bigg(\frac{1}{2}-\frac{1}{p_{2}},1-\alpha
_{1}-\alpha_{2}\bigg).$

\vspace{0.5cm} For the term $B_{1,3}^{1}(n,c)(t)$, using that $p_{3}\in
\lbrack4/3,2)$ and $\frac{1}{p_{1}}+\frac{1}{p_{3}}\leq\frac{3}{2}$ (see
(\ref{A2})), we can choose $q,q_{3}\in\lbrack1,\infty)$ such that $\frac
{1}{q_{3}}=\frac{1}{p_{3}}-\frac{1}{2}$ and $\frac{1}{q}=\frac{1}{p_{1}}%
+\frac{1}{q_{3}},$ and then proceed as follows%

\begin{align}
||B_{1,3}^{1}(n,\zeta)(t)||_{p_{1}}  &  =\bigg|\bigg|\int_{0}^{t}\nabla\cdot
e^{(t-s)\Delta}(n(S\ast\zeta))(s)ds\bigg|\bigg|_{p_{1}}\nonumber\\
&  \leq\int_{0}^{t}||\nabla\cdot e^{(t-s)\Delta}(n(S\ast\zeta))(s)||_{p_{1}%
}ds\nonumber\\
&  \leq C\int_{0}^{t}{(t-s)^{-\frac{1}{2}-(\frac{1}{q}-\frac{1}{p_{1}})}%
}||(n(S\ast\zeta))(s)||_{q}ds\nonumber\\
&  \leq C\int_{0}^{t}{(t-s)^{-\frac{1}{p_{3}}}}||n(s)||_{p_{1}}||S\ast
\zeta(s)||_{q_{3}}ds,\ \ \label{aux-B113-1}%
\end{align}
where above we have used estimate (\ref{est-sg-1}) and H\"{o}lder inequality.
Now, Hardy-Littlewood-Sobolev inequality (\ref{littlewood}), Beta functions
properties and the conditions $p_{3}\in\lbrack4/3,2)$ and (\ref{A4}) lead us
to
\begin{align}
\ \text{R.H.S. of (\ref{aux-B113-1})}  &  \leq C\sigma_{p_{3}}\int_{0}%
^{t}{(t-s)^{-\frac{1}{p_{3}}}}||n(s)||_{p_{1}}||\zeta(s)||_{p_{3}%
}ds\ \nonumber\\
&  \leq C\sigma_{p_{3}}|||n|||_{p_{1},\alpha_{1}}|||\zeta|||_{p_{3},\alpha
_{3}}\int_{0}^{t}{(t-s)^{-\frac{1}{p_{3}}}}s^{-\alpha_{1}-\alpha_{3}%
}ds\nonumber\\
&  \leq C\sigma_{p_{3}}|||n|||_{p_{1},\alpha_{1}}|||\zeta|||_{p_{3},\alpha
_{3}}t^{-\alpha_{1}}\beta\bigg(1-\frac{1}{p_{3}},1-\alpha_{1}-\alpha
_{3}\bigg)\ \nonumber\\
&  =C\sigma_{p_{3}}||n||_{X_{1}}||\zeta||_{X_{3}}t^{-\alpha_{1}}%
\beta\bigg(1-\frac{1}{p_{3}},1-\alpha_{1}-\alpha_{3}\bigg). \label{aux-B113-2}%
\end{align}
Putting together (\ref{aux-B113-1}) and (\ref{aux-B113-2}), considering
$C_{1,3}^{1}:=C\sigma_{p_{3}}\beta\bigg(1-\frac{1}{p_{3}},1-\alpha_{1}%
-\alpha_{3}\bigg),$ multiplying the resulting estimate by $t^{\alpha_{1}},$
and taking the supremum in $t$, we obtain the estimate
\[
||B_{1,3}^{1}(n,\zeta)||_{X_{1}}\leq C_{1,3}^{1}||n||_{X_{1}}||\zeta||_{X_{3}%
}.
\]

\vspace{0.5cm} Next we turn to the term $B_{2,3}^{2}(c,\zeta)(t)$. Here, using
$p_{3}\in\lbrack4/3,2)$ and $\frac{1}{p_{2}}+\frac{1}{p_{3}}<\frac{3}{2}$ (see
(\ref{A3})) and considering $\frac{1}{q_{3}}=\frac{1}{p_{3}}-\frac{1}{2}$ and
$\frac{1}{q}=\frac{1}{p_{2}}+\frac{1}{q_{3}},$ we can apply (\ref{est-sg-1}),
H\"{o}lder inequality and estimate (\ref{littlewood}) in order to obtain
\begin{align*}
||B_{2,3}^{2}(c,\zeta)(t)||_{\infty}  &  \leq\int_{0}^{t}||e^{(t-s)\Delta
}((S\ast\zeta)\cdot\nabla c)(s)||_{\infty}ds\\
&  \leq C\int_{0}^{t}{(t-s)^{-\frac{1}{q}}}||((S\ast\zeta)\cdot\nabla
c)(s)||_{q}ds,\ \ \\
&  \leq C\int_{0}^{t}{(t-s)^{-\frac{1}{p_{2}}-\frac{1}{p_{3}}+\frac{1}{2}}%
}||S\ast\zeta(s)||_{q_{3}}||\nabla c(s)||_{p_{2}}ds\ \ \ \\
&  \leq C\sigma_{p_{3}}\int_{0}^{t}{(t-s)^{-\frac{1}{p_{2}}-\frac{1}{p_{3}%
}+\frac{1}{2}}}||\zeta(s)||_{p_{3}}||\nabla c(s)||_{p_{2}}ds\ \ \ \\
&  \leq C\sigma_{p_{3}}|||\nabla c|||_{p_{2},\alpha_{2}}|||\zeta
|||_{p_{3},\alpha_{3}}\int_{0}^{t}{(t-s)^{-\frac{1}{p_{2}}-\frac{1}{p_{3}%
}+\frac{1}{2}}}s^{-\alpha_{2}-\alpha_{3}}ds\\
&  \leq C\sigma_{p_{3}}||c||_{X_{2}}||\zeta||_{X_{3}}t^{\frac{3}{2}-\frac
{1}{p_{2}}-\frac{1}{p_{3}}-\alpha_{2}-\alpha_{3}}\beta\bigg(\frac{3}{2}%
-\frac{1}{p_{2}}-\frac{1}{p_{3}},1-\alpha_{2}-\alpha_{3}\bigg),
\end{align*}
which implies that
\begin{equation}
|||B_{2,3}^{2}(c,\zeta)|||_{\infty,0}\leq\tilde{C_{1}}||c(s)||_{X_{2}}%
||\zeta||_{X_{3}}, \label{B223}%
\end{equation}
where $\tilde{C}_{1}:=C\sigma_{p_{3}}\mathcal{B}\bigg(\frac{3}{2}-\frac
{1}{p_{2}}-\frac{1}{p_{3}},1-\alpha_{2}-\alpha_{3}\bigg),$ because $\frac
{3}{2}-\frac{1}{p_{2}}-\frac{1}{p_{3}}-\alpha_{2}-\alpha_{3}=0$ by (\ref{A1}).

\vspace{0.5cm} For $\nabla B_{2,3}^{2}(c,\zeta)(t)$, again using $\frac{1}%
{q}=\frac{1}{p_{2}}+\frac{1}{q_{3}}$ and $\frac{1}{q_{3}}=\frac{1}{p_{3}%
}-\frac{1}{2},$ we can estimate%
\begin{align*}
||\nabla B_{2,3}^{2}(c,\zeta)(t)||_{p_{2}}  &  =\bigg|\bigg|\int_{0}^{t}\nabla
e^{(t-s)\Delta}((S\ast\zeta)\cdot\nabla c)(s)ds\bigg|\bigg|_{p_{2}}\\
&  \leq\int_{0}^{t}||\nabla e^{(t-s)\Delta}((S\ast\zeta)\cdot\nabla
c)(s)||_{p_{2}}ds\\
&  \leq C\int_{0}^{t}{(t-s)^{-\frac{1}{2}-\frac{1}{q}+\frac{1}{p_{2}}}%
}||((S\ast\zeta)\cdot\nabla c)(s)||_{q}ds\ \ \\
&  \leq C\int_{0}^{t}{(t-s)^{-\frac{1}{p_{3}}}}||S\ast\zeta(s)||_{q_{3}%
}||\nabla c(s)||_{p_{2}}ds\ \ \ \\
&  \leq C\sigma_{p_{3}}\int_{0}^{t}{(t-s)^{-\frac{1}{p_{3}}}}||\zeta
(s)||_{p_{3}}||\nabla c(s)||_{p_{2}}ds\ \ \ \ \text{by\ \ (\ref{littlewood}%
)}\\
&  \leq C\sigma_{p_{3}}|||\nabla c|||_{p_{2},\alpha_{2}}|||\zeta
|||_{p_{3},\alpha_{3}}\int_{0}^{t}{(t-s)^{-\frac{1}{p_{3}}}}s^{-\alpha
_{2}-\alpha_{3}}ds\\
&  \leq C\sigma_{p_{3}}||c||_{X_{2}}||\zeta||_{X_{3}}t^{-\alpha_{2}}%
\beta\bigg(1-\frac{1}{p_{3}},1-\alpha_{2}-\alpha_{3}\bigg),
\end{align*}
where in the last passage we use $\alpha_{3}=1-\frac{1}{p_{3}}$ (see
(\ref{A1})). Now, multiplying both sides of the inequality by $t^{\alpha_{2}}$
and taking the supremum over $t$, the resulting estimate is
\begin{equation}
|||\nabla B_{2,3}^{2}(c,\zeta)|||_{p_{2},\alpha_{2}}\leq\tilde{C_{2}%
}||c||_{X_{2}}||\zeta||_{X_{3}}, \label{GB223}%
\end{equation}
where $\tilde{C}:=C\sigma_{p_{3}}\beta\bigg(1-\frac{1}{p_{3}},1-\alpha
_{2}-\alpha_{3}\bigg).$ By (\ref{B223}) and (\ref{GB223}), and taking
$C_{2,3}^{2}:=\tilde{C_{1}}+\tilde{C_{2}},$ we get
\[
||B_{2,3}^{2}(c,\zeta)||_{X_{2}}\leq C_{2,3}^{2}||c||_{X_{2}}||\zeta||_{X_{3}%
}.
\]

\vspace{0.5cm} For the term $B_{1,2}^{2}(n,c)(t)$, using $\alpha_{1}%
=1-\frac{1}{p_{1}}$ (see (\ref{A1})) and proceeding similarly to estimate
(\ref{B223}) for $B_{2,3}^{2},$ it follows that
\begin{align*}
||B_{1,2}^{2}(n,c)(t)||_{\infty}  &  \leq C\int_{0}^{t}{(t-s)^{-\frac{1}%
{p_{1}}}}||n(s)||_{p_{1}}||c(s)||_{\infty}ds\ \ \ \\
&  \leq C|||n|||_{p_{1},\alpha_{1}}|||c|||_{\infty,0}\int_{0}^{t}%
{(t-s)^{-\frac{1}{p_{1}}}}s^{-\alpha_{1}}ds\\
&  \leq C||n||_{X_{1}}||c||_{X_{2}}\beta\bigg(1-\frac{1}{p_{1}},1-\alpha
_{1}\bigg)
\end{align*}
and then, considering $\tilde{C_{1}}:=C\beta\bigg(1-\frac{1}{p_{1}}%
,1-\alpha_{1}\bigg)$,
\begin{equation}
|||B_{1,2}^{2}(n,c)|||_{\infty,0}\leq\tilde{C_{1}}||n||_{X_{1}}||c||_{X_{2}}.
\label{B212}%
\end{equation}

\vspace{0.5cm}For the term $\nabla B_{1,2}^{2}(n,c)(t)$, we can proceed as in
the estimate for $\nabla B_{2,3}^{2}$ (see (\ref{GB223})) to obtain%
\begin{align*}
||\nabla B_{1,2}^{2}(n,c)(t)||_{p_{2}}  &  \leq C\int_{0}^{t}{(t-s)^{-\frac
{1}{2}-\frac{1}{p_{1}}+\frac{1}{p_{2}}}}||n(s)||_{p_{1}}||c(s)||_{\infty
}ds\ \ \\
&  \leq C|||n|||_{p_{1},\alpha_{1}}|||c|||_{\infty,0}\int_{0}^{t}%
{(t-s)^{-\frac{1}{2}-\frac{1}{p_{1}}+\frac{1}{p_{2}}}}s^{-\alpha_{1}}ds\\
&  \leq C||n||_{X_{1}}||c||_{X_{2}}t^{-\alpha_{2}}\beta\bigg(\frac{1}{2}%
-\frac{1}{p_{1}}+\frac{1}{p_{2}},1-\alpha_{1}\bigg),
\end{align*}
which yields
\begin{equation}
|||\nabla B_{1,2}^{2}(n,c)|||_{p_{2},\alpha_{2}}\leq\tilde{C_{2}}||n||_{X_{1}%
}||c||_{X_{2}}, \label{GB212}%
\end{equation}
where $\tilde{C_{2}}:=C\beta\bigg(\frac{1}{2}-\frac{1}{p_{1}}+\frac{1}{p_{2}%
},1-\alpha_{1}\bigg),$ after using the relations in (\ref{A1}), multiplying
the resulting inequality by $t^{\alpha_{2}}$, and computing the supremum over
$t.$ Now, estimates (\ref{B212}) and (\ref{GB212}) imply that
\[
|||B_{1,2}^{2}(c,\zeta)|||_{X_{2}}\leq C_{1,2}^{2}||c||_{X_{2}}||\zeta
||_{X_{3}},
\]
where $C_{1,2}^{2}:=\tilde{C_{1}}+\tilde{C_{2}}.$

In the sequel we treat with the bilinear operator $B_{3,3}^{3}.$ For that,
using that $p_{3}\in\lbrack4/3,2)$ we can take $\frac{1}{q_{3}}=\frac{1}%
{p_{3}}-\frac{1}{2}$ and $\frac{1}{q}=\frac{1}{p_{3}}+\frac{1}{q_{3}}$.
Passing the $L^{p_{3}}$-norm inside the integral, using (\ref{est-sg-1}) and
after H\"{o}lder inequality, we get%

\begin{align}
||B_{3,3}^{3}(\zeta,\tilde{\zeta})(t)||_{p_{3}}  &  \leq\int_{0}^{t}%
||\nabla\cdot e^{(t-s)\Delta}(\zeta(S\ast\tilde{\zeta}))(s)||_{p_{3}%
}ds\nonumber\\
&  \leq C\int_{0}^{t}{(t-s)^{-\frac{1}{2}-(\frac{1}{q}-\frac{1}{p_{3}})}%
}||(\zeta(S\ast\tilde{\zeta}))(s)||_{q}ds,\ \nonumber\\
&  \leq C\int_{0}^{t}{(t-s)^{-\frac{1}{p_{3}}}}||\zeta(s)||_{p_{3}}%
||S\ast\tilde{\zeta}(s)||_{q_{3}}ds.\ \ \ \ \label{aux-B333-1}%
\end{align}
Employing now (\ref{littlewood}) in (\ref{aux-B333-1}), recalling (\ref{beta})
and (\ref{A1}), we can estimate $B_{3,3}^{3}$ as follows%
\begin{align}
||B_{3,3}^{3}(\zeta,\tilde{\zeta})(t)||_{p_{3}}  &  \leq C\sigma_{p_{3}}%
\int_{0}^{t}{(t-s)^{-\frac{1}{p_{3}}}}||\zeta(s)||_{p_{3}}||\tilde{\zeta
}(s)||_{p_{3}}ds\nonumber\\
&  \leq C\sigma_{p_{3}}|||\zeta|||_{p_{3},\alpha_{3}}|||\tilde{\zeta
}|||_{p_{3},\alpha_{3}}\int_{0}^{t}{(t-s)^{-\frac{1}{p_{3}}}}s^{-2\alpha_{3}%
}ds\nonumber\\
&  =C\sigma_{p_{3}}||\zeta||_{X_{3}}||\tilde{\zeta}||_{X_{3}}t^{-\alpha_{3}%
}\beta\bigg(1-\frac{1}{p_{3}},1-2\alpha_{3}\bigg). \label{aux-B333-2}%
\end{align}
Using $\alpha_{3}<\frac{1}{2}$ and $p_{3}\in\lbrack4/3,2)$, taking
$C_{3,3}^{3}:=C\sigma_{p_{3}}\beta\bigg(1-\frac{1}{p_{3}},1-2\alpha_{3}%
\bigg)$, multiplying both sides of (\ref{aux-B333-2}) by $t^{\alpha_{3}},$ and
afterwards computing the supremum over $t$, we arrive at
\[
||B_{3,3}^{3}(\zeta,\tilde{\zeta})||_{X_{3}}\leq C_{3,3}^{3}||\zeta||_{X_{3}%
}||\tilde{\zeta}||_{X_{3}}.
\]

Finally, we handle the linear operator $L_{1}^{3}.$ Using that $p_{1}\geq2$,
we can take $q\in\lbrack1,\infty)$ such that $\frac{1}{q}=\frac{1}{p_{1}%
}+\frac{1}{2}$. Recalling that $0<\frac{1}{p_{3}}-\frac{1}{p_{1}}\leq\frac
{1}{2}$ (by (\ref{A3})), we can use (\ref{est-sg-1}), H\"{o}lder inequality
and (\ref{A1}) in order to estimate%

\begin{align*}
||L_{1}^{3}(n)(t)||_{p_{3}}  &  =\bigg|\bigg|\int_{0}^{t}\nabla^{\perp}\cdot
e^{(t-s)\Delta}n(s)\nabla\phi ds\bigg|\bigg|_{p_{3}}\\
&  \leq\int_{0}^{t}||\nabla^{\perp}\cdot e^{(t-s)\Delta}n(s)\nabla
\phi||_{p_{3}}ds\\
&  \leq C\int_{0}^{t}{(t-s)^{-\frac{1}{2}-(\frac{1}{q}-\frac{1}{p_{3}})}%
}||n(s)\nabla\phi||_{q}ds\ \ \ \\
&  \leq C\int_{0}^{t}{(t-s)^{-{1+}\frac{1}{p_{3}}-\frac{1}{p_{1}}}%
}||n(s)||_{p_{1}}||\nabla\phi||_{2}ds\ \ \\
&  \leq C|||n|||_{p_{1},\alpha_{1}}||\nabla\phi||_{2}\int_{0}^{t}%
{(t-s)^{-{1+}\frac{1}{p_{3}}-\frac{1}{p_{1}}}}s^{-\alpha_{1}}ds\\
&  =C||n||_{X_{1}}||\nabla\phi||_{2}t^{-\alpha_{3}}\beta\bigg(\frac{1}{p_{3}%
}-\frac{1}{p_{1}},1-\alpha_{1}\bigg),\
\end{align*}
which, taking $\alpha:=C||\nabla\phi||_{2}\beta\bigg(\frac{1}{p_{3}}-\frac
{1}{p_{1}},1-\alpha_{1}\bigg)$ and recalling $\alpha_{1}<1$ (by (\ref{A4})),
leads us to
\[
||L_{1}^{3}(n)||_{X_{3}}\leq\alpha||n||_{X_{1}}.
\]
In summary, the $B_{i,j}^{k}:X_{i}\times X_{j}\rightarrow X_{k}$ in
(\ref{IFN}) are all continuous bilinear maps with the constants $C_{i,j}^{k}$
of the corresponding bilinear estimate independent of $T$. Moreover,
$L_{1}^{3}:X_{1}\rightarrow X_{3}$ is a continuous linear map with constant
$\alpha$ independent of $T.$

In what follows we perform a fixed-point argument via Lemma \ref{LA}. \ For
that, we take $K_{1}:=1+\alpha$ and $K_{2}:=\alpha\sum_{i,j=1}^{3}C_{i,j}%
^{1}+\sum_{k,i,j=1}^{3}C_{i,j}^{k}$. Fix $\epsilon\in\big(0,\frac{1}%
{4K_{1}K_{2}}\big)$ and define $\mathcal{B}_{\epsilon}:=\{g\in X:||g||_{X}%
\leq2K_{1}\epsilon\}.$ Let $y$ be given by%
\[
y:=(y_{1},y_{2},y_{3})=(e^{t\Delta}n_{0},e^{t\Delta}c_{0},e^{t\Delta}\zeta
_{0}).
\]
Note that, by heat semigroup properties in Section \ref{sectionP} (see
(\ref{est-sg-1}) and (\ref{Prop-point-1})), we have that
\begin{equation}
e^{t\Delta}n_{0}\in\dot{C}_{\alpha_{1}}(L^{p_{1}}(\mathbb{R}^{2}%
)),\ \ e^{t\Delta}\nabla c_{0}\in\dot{C}_{\alpha_{2}}(L^{p_{2}}(\mathbb{R}%
^{2})),\ \ e^{t\Delta}\zeta_{0}\in\dot{C}_{\alpha_{3}}(L^{p_{3}}%
(\mathbb{R}^{2})). \label{AuxDotc}%
\end{equation}
By definition of $\dot{C}_{\alpha}$-spaces, there exists a $T_{0}>0$
(depending on $\epsilon$) such that
\begin{equation}
|||e^{t\Delta}n_{0}|||_{p_{1},\alpha_{1}}\leq\frac{\epsilon}{4}%
,\ \ |||e^{t\Delta}\nabla c_{0}|||_{p_{2},\alpha_{2}}\leq\frac{\epsilon}%
{4},\ \ |||e^{t\Delta}\zeta_{0}|||_{p_{3},\alpha_{3}}\leq\frac{\epsilon}{4}.
\label{aux-data-epsilon-1}%
\end{equation}
Recall the $X_{2}$-norm in (\ref{norm}), using (\ref{aux-data-epsilon-1}), and
the fact that $||c_{0}||_{\infty}$ is sufficiently small, then there exists a
$T_{0}>0$ such that (in $(0,T_{0})$)
\[
||y||_{X}=||(y_{1},y_{2},y_{3})||_{X}=||e^{t\Delta}n_{0}||_{X_{1}%
}+||e^{t\Delta}c_{0}||_{X_{2}}+||e^{t\Delta}\zeta_{0}||_{X_{3}}\leq\epsilon.
\]

Therefore, applying Lemma \ref{LA}, we obtain that the map $\Phi(n,c,\zeta)$
given in (\ref{Op-fi}) has a unique fixed-point $(n,c,\zeta)\in X$ such that
$||(n,c,\zeta)||_{X}\leq2K_{1}\epsilon,$ which is the mild solution for system
(\ref{CNS}), locally in $(0,T_{0}),$ with initial data $(n_{0},c_{0},\zeta
_{0})$.\hfill\hfill$\diamond$\vskip12pt

Some further comments on Theorem \ref{TLS} are in order.

\begin{remark}
\label{RTL}

\begin{itemize}
\item[$(i)$] (Large $L^{\infty}$-norm of $c_{0}$) For $s\in\mathbb{R}$ and
$1\leq p\leq\infty,$ let $\dot{H}_{p}^{s}(\mathbb{R}^{2})=\{f\in
\mathcal{S}^{\prime}/\mathcal{P};$ $\left\Vert f\right\Vert _{\dot{H}_{p}^{s}%
}=\left\Vert (-\Delta)^{s/2}f\right\Vert _{p}<\infty\}$ stand for the
homogeneous Sobolev space (see, e.g., \cite{Grafakos}), where $\mathcal{P}$
\ is the set of all polynomials in $\mathbb{R}^{2}$. In Theorem \ref{TLS}, we
could drop the condition $c_{0}\in L^{\infty}$ (and small $||c_{0}||_{\infty}%
$) and consider only $c_{0}\in\dot{H}_{2}^{1}(\mathbb{R}^{2})$, by modifying
the component space $X_{2}$ in such a way that
\[
X_{2}=\{c\in C((0,T);\dot{H}_{^{p_{2}}}^{1}(\mathbb{R}^{2}):\nabla
c\in{C_{\alpha_{2}}((0,T);}L^{p_{2}}(\mathbb{R}^{2}))\text{ and }%
(-\Delta)^{1/4}c\in C_{\widetilde{\alpha}_{2}}((0,T);L^{\widetilde{p}_{2}%
}(\mathbb{R}^{2}))\},
\]
endowed with the norm
\begin{equation}
||c||_{X_{2}}:=|||\nabla c|||_{p_{2},\alpha_{2}}+|||(-\Delta)^{1/4}%
c|||_{\widetilde{p}_{2},\widetilde{\alpha}_{2}}=\sup_{0<t<T}t^{\alpha_{2}%
}||\nabla c(t)||_{p_{2}}+\sup_{0<t<T}t^{\widetilde{\alpha}_{2}}\left\Vert
(-\Delta)^{1/4}c\right\Vert _{\widetilde{p}_{2}}, \label{aux-rem-space-1}%
\end{equation}
and adapting the corresponding estimates for the bilinear terms $B_{2,3}^{2}$
and $B_{1,2}^{2},$ where $4<\widetilde{p}_{2}<\infty$ and $\ \widetilde
{\alpha}_{2}+\frac{1}{\widetilde{p}_{2}}=\frac{1}{4}.$ For example, for
$B_{1,2}^{2},$ using the second parcel in (\ref{aux-rem-space-1}), we can
handle the product $\left\Vert nc\right\Vert _{q}\leq\left\Vert n\right\Vert
_{p_{1}}\left\Vert c\right\Vert _{r}$ and $\left\Vert c\right\Vert _{r}\leq
C\left\Vert (-\Delta)^{1/4}c\right\Vert _{\widetilde{p}_{2}},$ avoiding the
$L^{\infty}$-norm of $c$ and problems with the Sobolev embedding in the
critical case. For the linear part, similarly to (\ref{AuxDotc}%
)-(\ref{aux-data-epsilon-1}) and using that $\dot{H}_{2}^{1}(\mathbb{R}%
^{2})\hookrightarrow\dot{H}_{4}^{1/2}(\mathbb{R}^{2})$, we have that
\[
|||\nabla e^{t\Delta}c_{0}|||_{p_{2},\alpha_{2}},|||(-\Delta)^{1/4}e^{t\Delta
}c_{0}|||_{\widetilde{p}_{2},\widetilde{\alpha}_{2}}\leq C\left\Vert
c_{0}\right\Vert _{\dot{H}_{2}^{1}}%
\]
as well as they can be taken sufficiently small when $T$ is small. Moreover,
considering an arbitrary $c_{0}\in\dot{H}_{2}^{1}\cap L^{\infty},$ then we
obtain additionally that $c\in C_{0}((0,T);L^{\infty}(\mathbb{R}^{2})).$
However, the proof would be somewhat more elaborate. Anyway, in order to
extend the solution globally, we will need a smallness condition on
$||c_{0}||_{\infty}$. Thus, for simplicity, we prefer to already present the
local-in-time results with that condition on $||c_{0}||_{\infty}.$

\item[$(ii)$] Another way to drop the smallness condition on $||c_{0}%
||_{\infty}$ would be work with non-critical spaces for $n$ and $c$. In fact,
for $\delta\in(0,1)$, considering the component spaces $X_{1},X_{2}$ in
(\ref{spacesolutions-0})-(\ref{spacesolutions}) with $\alpha_{1}+\frac
{1}{p_{1}}=1-\delta$ and
\[
X_{2}=\{c\in C_{\delta}((0,T);L^{\infty}(\mathbb{R}^{2})),\text{ }\nabla
c\in{C_{\alpha_{2}}((0,T);}L^{p_{2}}(\mathbb{R}^{2})\},
\]
it can be shown that there exists a $T_{0}>0$ such that (\ref{CNS}) has a
unique mild solution $(n,c,\zeta)\in X$, locally in $(0,T_{0})$, with initial
data $n_{0}\in L^{\frac{1}{1-\delta}}(\mathbb{R}^{2}),$ $c_{0}\in L^{\infty
}(\mathbb{R}^{2})$ with $\nabla c_{0}\in L^{2}(\mathbb{R}^{2})$ and $\zeta
_{0}\in L^{1}(\mathbb{R}^{2})$ or, according to Section \ref{sectionM},
$\zeta_{0}\in\mathcal{M}(\mathbb{R}^{2})$.

\item[$(iii)$] (General parameter functions) The proof of Theorem \ref{TLS}
and corresponding estimates could be adapted to consider the more general
cases of chemotactic sensitivity and oxygen consumption rate. In fact, we can
consider $\chi,f\in C^{1}(\mathbb{R}^{+})$ satisfying (\ref{Basic-cond-1}).
For example, in order to handle the product $nf(c)$ in (\ref{CNS}) in the
contraction argument, we should use $\left\vert f(c)\right\vert \leq M|c|$ and
$\left\vert f(c)-f(\tilde{c})\right\vert \leq M|c-\tilde{c}|$ which is valid
with $M=f^{\prime}(\Gamma)$ when $|||c|||_{\infty,0}\leq\Gamma.$
\end{itemize}
\end{remark}

In the next proposition we extend the range of indexes $p_{1},p_{2},p_{3}$ for
which the solution belongs to class (\ref{spacesolutions-0}%
)-(\ref{spacesolutions}).

\begin{proposition}
\label{PropLp} (Further $L^{p}$-integrability) Let $(n,c,\zeta)$ be the
solution obtained in Theorem \ref{TLS}. Then, $(n,c,\zeta)$ belongs to the
class (\ref{spacesolutions-0})-(\ref{spacesolutions}), for every $q_{1}%
,q_{3}\in\lbrack1,\infty)$, $q_{2}\in\lbrack2,\infty)$ and $\alpha_{1}%
,\alpha_{2},\alpha_{3}\geq0$ such that
\begin{equation}
\alpha_{1}+\frac{1}{q_{1}}=1,\ \ \alpha_{2}+\frac{1}{q_{2}}=\frac{1}%
{2},\ \ \alpha_{3}+\frac{1}{q_{3}}=1. \label{aux-rem-iii}%
\end{equation}

\end{proposition}

\textit{Proof.} Using the integral formulation (\ref{IFN}), we can write $n$
as
\begin{equation}
n(t)=e^{t\Delta}n_{0}-B_{1,2}^{1}(n,c)(t)-B_{1,3}^{1}(n,\zeta)(t).
\label{aux-n-proof-1}%
\end{equation}
Firstly, let $\alpha_{1},\alpha_{2},\alpha_{3},p_{1},p_{2},p_{3}$ be as in
\textbf{Condition A}. We are going to show that $n\in C_{\alpha}(L^{q})$, for
all $q\geq1$ and $\alpha\in\lbrack0,1)$ such that
\begin{equation}
\alpha+\frac{1}{q}=1. \label{PropLp01}%
\end{equation}

By (\ref{est-sg-1}) we have that $e^{t\Delta}n_{0}\in C_{\alpha}(L^{q})$, for
all $q,\alpha$ as above. For the second term in the R.H.S of
(\ref{aux-n-proof-1}), let $\frac{1}{r_{1}}=\frac{1}{p_{1}}+\frac{1}{p_{2}}$.
Then, for $r_{1}\leq q$, it follows that%

\begin{align}
||B_{1,2}^{1}(n,c)(t)||_{q}  &  \leq\int_{0}^{t}||\nabla\cdot e^{(t-s)\Delta
}(n\nabla c)(s)||_{q}ds\nonumber\\
&  \leq C\int_{0}^{t}{(t-s)^{-\frac{1}{2}-(\frac{1}{r_{1}}-\frac{1}{q})}%
}||(n\nabla c)(s)||_{r_{1}}ds\nonumber\\
&  \leq C\int_{0}^{t}{(t-s)^{{\frac{1}{q}-\frac{1}{2}-\frac{1}{p_{1}}-\frac
{1}{p_{2}}}}}||n(s)||_{p_{1}}||\nabla c(s)||_{p_{2}}ds\nonumber\\
&  \leq C|||n|||_{p_{2},\alpha_{1}}|||\nabla c|||_{p_{2},\alpha_{2}}\int
_{0}^{t}{(t-s)^{{\frac{1}{q}-\frac{1}{2}-\frac{1}{p_{1}}-\frac{1}{p_{2}}}}%
}s^{-\alpha_{1}-\alpha_{2}}ds.\ \ \ \ \label{REaux-B112-1}%
\end{align}
Assuming now that ${-\frac{1}{q}+\frac{1}{2}+\frac{1}{p_{1}}+\frac{1}{p_{2}}%
}<1$ \big(i.e., $q<\frac{1}{1-\alpha_{1}-\alpha_{2}}$\big), we can use
(\ref{beta}) and \textbf{Condition A} to obtain
\begin{align}
||B_{1,2}^{1}(n,c)(t)||_{q}  &  \leq C|||n|||_{p_{2},\alpha_{1}}|||\nabla
c|||_{p_{2},\alpha_{2}}t^{\frac{1}{2}+\frac{1}{q}-\frac{1}{p_{1}}-\frac
{1}{p_{2}}-\alpha_{1}-\alpha_{2}}\beta\bigg(\frac{1}{q}+\frac{1}{2}-\frac
{1}{p_{1}}-\frac{1}{p_{2}},1-\alpha_{1}-\alpha_{2}\bigg)\ \ \ \ \nonumber\\
&  \leq C|||n|||_{p_{1},\alpha_{1}}|||\nabla c|||_{p_{2},\alpha_{2}}%
t^{\frac{1}{q}-1}\beta\bigg(\frac{1}{q}+\frac{1}{2}-\frac{1}{p_{1}}-\frac
{1}{p_{2}},1-\alpha_{1}-\alpha_{2}\bigg)\ \ \text{by}\ \text{\ (\ref{A1}%
)}\nonumber\\
&  \leq C||n||_{X_{1}}||c||_{X_{2}}t^{\frac{1}{q}-1}\beta\bigg(\frac{1}%
{q}+\frac{1}{2}-\frac{1}{p_{1}}-\frac{1}{p_{2}},1-\alpha_{1}-\alpha
_{2}\bigg)\text{.}\ \label{REaux-B112-2}%
\end{align}
Hence, we see that
\[
B_{1,2}^{1}(n,c)(t)\in C_{\alpha}(L^{q}),
\]
for all $q\geq1$ and $\alpha\in\lbrack0,1)$ satisfying (\ref{PropLp01}) and
such that $q\in\big[\frac{p_{1}p_{2}}{p_{1}+p_{2}},\frac{1}{1-\alpha
_{1}-\alpha_{2}}\big)$.

Finally, for the third term $B_{1,3}^{1}(n,\zeta)(t)$, we consider $\frac
{1}{r_{2}}=\frac{1}{p_{3}}-\frac{1}{2}$ and $\frac{1}{r_{3}}=\frac{1}{p_{1}%
}+\frac{1}{r_{2}}$. Then, for $r_{3}\leq q$, we can estimate%

\begin{align}
||B_{1,3}^{1}(n,\zeta)(t)||_{q}  &  \leq\int_{0}^{t}||\nabla\cdot
e^{(t-s)\Delta}(n(S\ast\zeta))(s)||_{q}ds\nonumber\\
&  \leq C\int_{0}^{t}{(t-s)^{-\frac{1}{2}-(\frac{1}{r_{3}}-\frac{1}{q})}%
}||(n(S\ast\zeta))(s)||_{r_{3}}ds\nonumber\\
&  \leq C\int_{0}^{t}{(t-s)^{\frac{1}{q}-\frac{1}{p_{1}}-\frac{1}{p_{3}}}%
}||n(s)||_{p_{1}}||S\ast\zeta(s)||_{r_{2}}ds.\ \ \label{REaux-B113-1}%
\end{align}
Assuming now that ${-\frac{1}{q}+\frac{1}{p_{1}}+\frac{1}{p_{3}}}<1$
\big(i.e., $q<\frac{1}{1-\alpha_{1}-\alpha_{3}}$\big), we can use
(\ref{beta}), (\ref{littlewood}) and \textbf{Condition A} to obtain
\begin{align}
||B_{1,3}^{1}(n,\zeta)(t)||_{q}  &  \leq C\sigma_{p_{3}}\int_{0}%
^{t}{(t-s)^{\frac{1}{q}-\frac{1}{p_{1}}-\frac{1}{p_{3}}}}||n(s)||_{p_{1}%
}||\zeta(s)||_{p_{3}}ds\ \nonumber\\
&  \leq C\sigma_{p_{3}}|||n|||_{p_{1},\alpha_{1}}|||\zeta|||_{p_{3},\alpha
_{3}}\int_{0}^{t}{(t-s)^{\frac{1}{q}-\frac{1}{p_{1}}-\frac{1}{p_{3}}}%
}s^{-\alpha_{1}-\alpha_{3}}ds\nonumber\\
&  \leq C\sigma_{p_{3}}|||n|||_{p_{1},\alpha_{1}}|||\zeta|||_{p_{3},\alpha
_{3}}t^{1+\frac{1}{q}-\frac{1}{p_{1}}-\frac{1}{p_{3}}-\alpha_{1}-\alpha_{3}%
}\beta\bigg(1+\frac{1}{q}-\frac{1}{p_{1}}-\frac{1}{p_{3}},1-\alpha_{1}%
-\alpha_{3}\bigg)\ \nonumber\\
&  =C\sigma_{p_{3}}||n||_{X_{1}}||\zeta||_{X_{3}}t^{\frac{1}{q}-1}%
\beta\bigg(1+\frac{1}{q}-\frac{1}{p_{1}}-\frac{1}{p_{3}},1-\alpha_{1}%
-\alpha_{3}\bigg). \label{REaux-B113-2}%
\end{align}
Noting that $r_{3}\leq q$ implies $q\geq\frac{2p_{1}p_{3}}{2p_{1}+2p_{3}%
-p_{1}p_{3}}$, we see that
\[
B_{1,3}^{1}(n,\zeta)(t)\in C_{\alpha}(L^{q}),
\]
for all $q\geq1$ and $\alpha\in\lbrack0,1)$ satisfying (\ref{PropLp01}) and
such that $q\in\big[\frac{2p_{1}p_{3}}{2p_{1}+2p_{3}-p_{1}p_{3}},\frac
{1}{1-\alpha_{1}-\alpha_{3}}\big)$.

Collecting the facts above, we arrive at $n\in C_{\alpha}(L^{q}),$ for all
$q\geq1$ and $\alpha\in\lbrack0,1)$ such that $\alpha+\frac{1}{q}=1$ and
\begin{equation}
q\in\bigg[\frac{p_{1}p_{2}}{p_{1}+p_{2}},\frac{1}{1-\alpha_{1}-\alpha_{2}%
}\bigg)\bigcap\bigg[\frac{2p_{1}p_{3}}{2p_{1}+2p_{3}-p_{1}p_{3}},\frac
{1}{1-\alpha_{1}-\alpha_{3}}\bigg). \label{PropLp02}%
\end{equation}
Proceeding similarly as done for the component $n$, we can also show that
$\nabla c\in C_{\alpha}(L^{q}),$ for all $q\geq2$ and $\alpha\in\lbrack0,1)$
such that $\alpha+\frac{1}{q}=\frac{1}{2}$ and
\begin{equation}
q\in\bigg[\frac{2p_{2}p_{3}}{2p_{2}+2p_{3}-p_{2}p_{3}},\frac{p_{2}p_{3}}%
{p_{2}+p_{3}-p_{2}p_{3}}\bigg)\bigcap\bigg[p_{1},\frac{2p_{1}}{2-p_{3}%
}\bigg) \label{PropLp02.1}%
\end{equation}
and that $\zeta\in C_{\alpha}(L^{q}),$ for all $q\geq1$ and $\alpha\in
\lbrack0,1)$ such that $\alpha+\frac{1}{q}=1$ and
\begin{equation}
q\in\bigg[\frac{2p_{1}}{2+p_{1}},p_{1}\bigg)\bigcap\bigg[\frac{2p_{3}}%
{4-p_{3}},\frac{p_{3}}{2-p_{3}}\bigg). \label{PropLp02.2}%
\end{equation}

Having established that $n\in C_{\alpha}(L^{q}),$ for all $q\geq1$ and
$\alpha\in\lbrack0,1)$ satisfying (\ref{PropLp02}) and $\alpha+\frac{1}{q}=1$,
we can repeat the above process finitely and then obtain the same property for
all $q\geq1$ and $\alpha\in\lbrack0,1)$ with $\alpha+\frac{1}{q}=1$. The same
idea can also be applied to $\nabla c$ and $\zeta$. For the reader
convenience, let us illustrate how this bootstrap-like process works. Consider
for example, the indexes $p_{1}=\frac{17}{8}$, $p_{2}=\frac{17}{8}$,
$p_{3}=\frac{15}{8}$, $\alpha_{1}=\frac{9}{17}$, $\alpha_{2}=\frac{1}{34}$ and
$\alpha_{3}=\frac{7}{15}$ which satisfy \textbf{Condition A}. Using
(\ref{PropLp02}), we obtain that $n\in C_{\alpha}(L^{q}),$ for all
\begin{equation}
q\in\lbrack289/272,34/15)\cap\lbrack510/257,255)=[510/257,34/15)
\label{PropLp03}%
\end{equation}
and $\alpha\in\lbrack0,1)$ such that $\alpha+\frac{1}{q}=1$. Considering
(\ref{PropLp02.1}), we have that $\nabla c\in C_{\alpha}(L^{q}),$ for all
\begin{equation}
q\in\lbrack510/257,255)\cap\lbrack17/8,34)=[17/8,34) \label{PropLp03.1}%
\end{equation}
and $\alpha\in\lbrack0,1)$ such that $\alpha+\frac{1}{q}=\frac{1}{2}$. In
turn, using (\ref{PropLp02.2}), it follows that $\zeta\in C_{\alpha}(L^{q}),$
for all
\begin{equation}
q\in\lbrack34/33,17/8)\cap\lbrack30/17,15)=[30/17,17/8) \label{PropLp03.2}%
\end{equation}
and $\alpha\in\lbrack0,1)$ such that $\alpha+\frac{1}{q}=1$.

Note that using the new ranges obtained in (\ref{PropLp03}), (\ref{PropLp03.1}%
) and (\ref{PropLp03.2}), we can repeat the same process by starting now with
$p_{1}=q_{1},$ $p_{2}=q_{2},$ and $p_{3}=q_{3}$ where $q_{1},q_{2},$ and
$q_{3}$ are values selected from the intervals $[510/257,34/15)$, $[17/8,34)$,
and $[30/17,17/8),$ respectively. For example, take $p_{1}=\frac{17}{8}$,
$p_{2}=17$ and $p_{3}=\frac{34}{19}$, then $\alpha_{1}=\frac{9}{17}$,
$\alpha_{2}=\frac{1}{34}$ and $\alpha_{3}=\frac{15}{34}$. So, applying
(\ref{PropLp02}), we obtain that $n\in C_{\alpha}(L^{q}),$ for all
\begin{equation}
q\in\lbrack289/153,34)\cap\lbrack289/153,34)=[289/153,34) \label{PropLp04}%
\end{equation}
and $\alpha\in\lbrack0,1)$ such that $\alpha+\frac{1}{q}=1$. Summarizing the
process up to this point (after only two iterations), starting from the value
$p_{1}=17/8$, we have progressed to the stage where $q$ falls in the interval
$[289/153,34)$. Now, one can iteratively apply this process a finite number of
times to ultimately achieve $n\in C_{\alpha}(L^{q})$ with $q\geq1$ and
$\alpha\in\lbrack0,1)$ such that $\alpha+\frac{1}{q}=1$. \hfill\hfill
$\diamond$\vskip12pt

\begin{remark}
\label{Rem-RTL-PropLp}

\begin{itemize}
\item[$(i)$] (Persistence property) In view of Proposition \ref{PropLp}, we
can consider $p_{1}=p_{3}=1$ and $\alpha_{1}=\alpha_{3}=0$ which leads us to
solutions $n\in BC([0,T_{0});L^{1}(\mathbb{R}^{2}))$ and $\xi\in
BC([0,T_{0});L^{1}(\mathbb{R}^{2}))$. Also, we already had that $c\in
BC((0,T_{0});L^{\infty}(\mathbb{R}^{2}))$. Assuming further that $c_{0}$
belongs to the space of bounded and uniformly continuous functions
$BUC(\mathbb{R}^{2})$, we get the strong continuity at $t=0^{+}$ for $c$ and
then $c\in BC([0,T_{0});L^{\infty}(\mathbb{R}^{2})).$

\item[$(ii)$] Theorem \ref{TLS} deals with $n_{0}\in L^{1}(\mathbb{R}^{2})$,
$c_{0}\in L^{\infty}(\mathbb{R}^{2})$, $\nabla c_{0}\in L^{2}(\mathbb{R}^{2})$
and $\zeta_{0}\in L^{1}(\mathbb{R}^{2})$, which are critical cases for the
initial-data spaces. With a suitable adaptation on the statement and proofs,
and keeping in mind Proposition \ref{PropLp}, a version of Theorem \ref{TLS}
is still valid in the subcritical case $n_{0}\in L^{p_{1}}(\mathbb{R}^{2})$,
$c_{0}\in L^{\infty}(\mathbb{R}^{2})$, $\nabla c_{0}\in L^{p_{2}}%
(\mathbb{R}^{2})$ and $\zeta_{0}\in L^{p_{3}}(\mathbb{R}^{2})$, for every
$p_{1},p_{3}>1$ and $p_{2}>2$ and $\alpha_{1},\alpha_{2},\alpha_{3}>0$
satisfying (\ref{aux-rem-iii}). To see this, we only have to notice that
$|||e^{t\Delta}n_{0}|||_{p_{1}}\leq||n_{0}||_{p_{1}}$, and then $|||e^{t\Delta
}n_{0}|||_{p_{1},\alpha_{1}}\leq T^{\alpha_{1}}||n_{0}||_{p_{1}}$ in $(0,T)$,
which can be made sufficiently small according to the size of $T>0$. The same
argument can be applied to $\nabla c_{0}\in L^{p_{2}}(\mathbb{R}^{2})$ and
$\zeta_{0}\in L^{p_{3}}(\mathbb{R}^{2})$.
\end{itemize}
\end{remark}

In what follows, exploiting a kind of parabolic regularization effect, we show
that the mild solutions obtained in Theorem \ref{TLS} are indeed smooth for
$t>0$.

\begin{proposition}
\label{regularity} (Regularity of solutions) Let $(n,c,\zeta)$ be the mild
solution of (\ref{CNS}) obtained in Theorem \ref{TLS}. Then, $(n,c,\zeta)$ is
smooth for $t>0$. More precisely, given $m_{1},m_{2}\in\mathbb{N}\cup\{0\}$,
if $\nabla\phi\in H^{m_{2}}(\mathbb{R}^{2})$ then
\begin{equation}
\partial_{t}^{m_{1}}\partial_{x}^{m_{2}}(n,c,\zeta)\in\widetilde{X}%
_{1}(0,T_{0}) \times\widetilde{X}_{2}(0,T_{0})\times\widetilde{X}_{3}(0,T_{0})
\label{Aux8.1}%
\end{equation}
where \text{ }
\[
\widetilde{X}_{1}(0,T_{0})=C((0,T_{0});L^{q_{1}}(\mathbb{R}^{2})),\text{
}\widetilde{X}_{2}(0,T_{0})=C((0,T_{0});L^{q_{2}}(\mathbb{R}^{2}))\text{ and
}\widetilde{X} _{3}(0,T_{0})=C((0,T_{0});L^{q_{3}}(\mathbb{R}^{2})),
\]
for every $q_{1},q_{3}\in\lbrack1,\infty)$ and $q_{2}\in\lbrack2,\infty)$.
\end{proposition}

\textit{Proof.} We are going to prove (\ref{Aux8.1}) by induction. First, we
consider the spatial derivatives, that is, the case $m_{1}=0$, and the indexes
$p_{1},p_{2}$ and $p_{3}$ satisfying \textbf{Condition A}.

By Theorem \ref{TLS}, we know that (\ref{Aux8.1}) holds for $m_{2}=0$.
Supposing (induction hypothesis) that (\ref{Aux8.1}) holds for $m_{2}=m-1$, we
are going to show it for $m_{2}=m$.

Since we are interested in times $t>0$, it is enough to take an arbitrary
$\sigma\in(0,T_{0})$ and show (\ref{Aux8.1}) in $(\sigma,T_{0})$. For
$\sigma,T\in(0,T_{0})$, with $\sigma<T$, define the sets
\[
X_{0}[\sigma,T):=C_{0}([\sigma,T);L^{p_{1}}(\mathbb{R}^{2}))\times\{c\in
C_{0}([\sigma,T);L^{\infty}(\mathbb{R}^{2}):\nabla c\in{C_{0}([\sigma
,T);L^{p_{2}}(\mathbb{R}^{2}))}\}\times C_{0}([\sigma,T);L^{p_{3}}%
(\mathbb{R}^{2}))
\]
and
\[
X_{\frac{1}{2}}(\sigma,T):=C_{\frac{1}{2}}((\sigma,T);L^{p_{1}}(\mathbb{R}%
^{2}))\times\{c\in C((\sigma,T);L^{\infty}(\mathbb{R}^{2}):\nabla
c\in{C_{\frac{1}{2}}((\sigma,T);L^{p_{2}}(\mathbb{R}^{2}))}\}\times
C_{\frac{1}{2}}((\sigma,T);L^{p_{3}}(\mathbb{R}^{2})).
\]

Given $K>0$, consider $B_{K}$ the subset of $X_{1}(\sigma,T)\times
X_{2}(\sigma,T)\times X_{3}(\sigma,T)$ consisting of all triples $(n,c,\zeta)$
such that $\partial_{x}^{l}(n,c,\zeta)\in X_{0}[\sigma,T)$, for all $0\leq
l\leq m-1,$ $\partial_{x}^{m}(n,c,\zeta)\in X_{\frac{1}{2}}(\sigma,T)$ and
\[
||\partial_{x}^{l}(n,c,\zeta)||_{X_{0}[\sigma,T)}\leq K,\text{ for}\ 0\leq
l\leq m-1,\text{ and}\ \ ||\partial_{x}^{m}(n,c,\zeta)||_{X_{\frac{1}{2}%
}(\sigma,T)}\leq K.
\]
By induction hypothesis, we have that $\partial_{x}^{l}(n(\sigma
),c(\sigma),\zeta(\sigma))\in L^{p_{1}}(\mathbb{R}^{2})\times L^{\infty
}(\mathbb{R}^{2})\cap\dot{H}^{1}(\mathbb{R}^{2})\times L^{p_{3}}%
(\mathbb{R}^{2})$ for all $0\leq l\leq m-1$. Then, it follows that
\begin{equation}
(e^{t\Delta}n(\sigma),e^{t\Delta}c(\sigma),e^{t\Delta}\zeta(\sigma))\in B_{K},
\label{AuxK}%
\end{equation}
for an appropriately chosen constant $K>0$.

In the sequel, for $\epsilon>0$ and $(n,c,\zeta)\in B_{K}$, we are going to
show that
\[
\label{Aux20}%
\begin{cases}
||\partial_{x}^{l}B_{1,2}^{1}(n,c)||_{C_{0}([\sigma,T);L^{p_{1}}%
(\mathbb{R}^{2}))}\leq\epsilon,\ \ \forall\ \ 0\leq l\leq m-1,\\
||\partial_{x}^{m}B_{1,2}^{1}(n,c)||_{C_{\frac{1}{2}}((\sigma,T);L^{p_{1}%
}(\mathbb{R}^{2}))}\leq\epsilon,
\end{cases}
\]
provided that $T-\sigma$ is small enough. Indeed, note initially that for
$l\in\mathbb{N}$, we have
\[
\partial_{x}^{l}B_{1,2}^{1}(n,c)(t)=\sum_{j=0}^{l}{\binom{l}{j}}B_{1,2}%
^{1}(\partial_{x}^{j}n,\partial_{x}^{l-j}c)(t).
\]
Then, for $0\leq l\leq m-1$, we can proceed as follows%

\begin{align*}
||\partial_{x}^{l}B_{1,2}^{1}(n,c)(t)||_{p_{1}}  &  \leq\sum_{j=0}^{l}%
{\binom{l}{j}}\int_{\sigma}^{t}||\nabla\cdot e^{(t-s)\Delta}(\partial_{x}%
^{j}n\partial_{x}^{l-j}\nabla c)(s)||_{p_{1}}ds\\
&  \leq C\sum_{j=0}^{l}{\binom{l}{j}}\int_{\sigma}^{t}{(t-s)^{-\frac{1}%
{2}-(\frac{1}{q}-\frac{1}{p_{1}})}}||(\partial_{x}^{j}n\partial_{x}%
^{l-j}\nabla c)(s)||_{q}ds,\ \ \text{with}\ \ \frac{1}{q}=\frac{1}{p_{1}%
}+\frac{1}{p_{2}}\\
&  \leq C\sum_{j=0}^{l}{\binom{l}{j}}\int_{\sigma}^{t}{(t-s)^{-\frac{1}%
{2}-\frac{1}{p_{1}}}}||\partial_{x}^{j}n(s)||_{p_{1}}||\partial_{x}%
^{l-j}\nabla c(s)||_{p_{2}}ds\\
&  \leq2^{l}CK^{2}\int_{\sigma}^{t}{(t-s)^{-\frac{1}{2}-\frac{1}{p_{1}}}}ds\\
&  \leq2^{l}CK^{2}\beta\bigg(\frac{1}{2}-\frac{1}{p_{1}},1\bigg)(t-\sigma
)^{\frac{1}{2}-\frac{1}{p_{1}}}.
\end{align*}
Now, calculating the supremum over $t\in\lbrack\sigma,T)$, we arrive at
\[
||\partial_{x}^{l}B_{1,2}^{1}(n,c)||_{C_{0}([\sigma,T);L^{p_{1}}%
(\mathbb{R}^{2}))}\leq2^{l}CK^{2}\beta\bigg(\frac{1}{2}-\frac{1}{p_{1}%
},1\bigg)(T-\sigma)^{\frac{1}{2}-\frac{1}{p_{1}}}.
\]
For $\partial_{x}^{m}B_{1,2}^{1}(n,c)(t)=\sum_{j=0}^{m}{\binom{m}{j}}%
B_{1,2}^{1}(\partial_{x}^{j}n,\partial_{x}^{m-j}c)(t)$, the terms $B_{1,2}%
^{1}(\partial_{x}^{j}n,\partial_{x}^{m-j}c)(t)$ with $0<j<m$ can be estimated
in the same way as in the previous case in order to obtain
\[
\bigg|\bigg|B_{1,2}^{1}(\partial_{x}^{j}n,\partial_{x}^{m-j}%
c)\bigg|\bigg|_{p_{1}}\leq CK^{2}\beta\bigg(\frac{1}{2}-\frac{1}{p_{1}%
},1\bigg)(t-\sigma)^{\frac{1}{2}-\frac{1}{p_{1}}}.
\]
Multiplying both sides of the above inequality by $(t-\sigma)^{\frac{1}{2}}$
and taking the supremum over $t\in\lbrack\sigma,T)$, we get%

\[
||B_{1,2}^{1}(\partial^{j}_{x}n,\partial^{m-j}_{x}c)||_{C_{\frac{1}{2}%
}((\sigma,T);L^{p_{1}}(\mathbb{R}^{2}))}\leq CK^{2}\beta\bigg(\frac{1}%
{2}-\frac{1}{p_{1}},1\bigg)(T-\sigma)^{1-\frac{1}{p_{1}}}.
\]
\newline For the case $j=m$, we have%

\begin{align*}
||B_{1,2}^{1}(\partial_{x}^{m}n,c)(t)||_{p_{1}}  &  \leq\int_{\sigma}%
^{t}||\nabla\cdot e^{(t-s)\Delta}(\partial_{x}^{m}n\nabla c)(s)||_{p_{1}}ds\\
&  \leq C\int_{\sigma}^{t}{(t-s)^{-\frac{1}{2}-(\frac{1}{q}-\frac{1}{p_{1}})}%
}||(\partial_{x}^{m}n\nabla c)(s)||_{q}ds\ \ \text{with }\frac{1}{q}=\frac
{1}{p_{1}}+\frac{1}{p_{2}}\\
&  \leq C\int_{\sigma}^{t}{(t-s)^{-\frac{1}{2}-\frac{1}{p_{1}}}}||\partial
_{x}^{m}n(s)||_{p_{1}}||\nabla c(s)||_{p_{2}}ds\\
&  \leq C\int_{\sigma}^{t}{(t-s)^{-\frac{1}{2}-\frac{1}{p_{1}}}}s^{-\frac
{1}{2}}s^{\frac{1}{2}}||\partial_{x}^{m}n(s)||_{p_{1}}||\nabla c(s)||_{p_{2}%
}ds\\
&  \leq CK^{2}\int_{\sigma}^{t}{(t-s)^{-\frac{1}{2}-\frac{1}{p_{1}}}}%
s^{-\frac{1}{2}}ds\\
&  \leq CK^{2}\beta\bigg(\frac{1}{2}-\frac{1}{p_{1}},\frac{1}{2}%
\bigg)(T-\sigma)^{-\frac{1}{p_{1}}},
\end{align*}
which yields
\[
||B_{1,2}^{1}(\partial_{x}^{m}n,c)||_{C_{\frac{1}{2}}((\sigma,T);L^{p_{1}%
}(\mathbb{R}^{2}))}\leq CK^{2}\beta\bigg(\frac{1}{2}-\frac{1}{p_{1}},\frac
{1}{2}\bigg)(T-\sigma)^{\frac{1}{2}-\frac{1}{p_{1}}}.
\]
Similarly, for $j=0$, we can show that%

\[
||B_{1,2}^{1}(n,\partial_{x}^{m}c)||_{C_{\frac{1}{2}}((\sigma,T);L^{p_{1}%
}(\mathbb{R}^{2}))}\leq CK^{2}\beta\bigg(\frac{1}{2}-\frac{1}{p_{1}},\frac
{1}{2}\bigg)(T-\sigma)^{\frac{1}{2}-\frac{1}{p_{1}}}.
\]
\newline So, taking $(T-\sigma)$ small enough leads us to%
\[
||\partial_{x}^{l}B_{1,2}^{1}(n,c)||_{C_{0}([\sigma,T);L^{p_{1}}%
(\mathbb{R}^{2}))}\leq\epsilon,\text{ }\forall\ \ 0\leq l\leq m-1,
\]
and
\[
||\partial_{x}^{m}B_{1,2}^{1}(n,c)||_{C_{\frac{1}{2}}((\sigma,T);L^{p_{1}%
}(\mathbb{R}^{2}))}\leq\epsilon.
\]
\newline Proceeding in the same way as above, we can show that%
\begin{equation}%
\begin{cases}
||\partial_{x}^{l}B_{i,j}^{k}(x_{i},x_{j})||_{C_{0}([\sigma,T);L^{p_{k}%
}(\mathbb{R}^{2}))}\leq\epsilon,\ \ \forall\ \ 0\leq l\leq m-1,\\
||\partial_{x}^{m}B_{i,j}^{k}(x_{i},x_{j})||_{C_{\frac{1}{2}}((\sigma
,T);L^{p_{k}}(\mathbb{R}^{2}))}\leq\epsilon,
\end{cases}
\label{Aux21}%
\end{equation}
where $(x_{1},x_{2},x_{3})=(n,c,\zeta)\in B_{K}$ and $(T-\sigma)$ is small
enough. \newline

Using that $\nabla\phi\in H^{m}(\mathbb{R}^{2})$ (which is true for all $m\leq
m_{2}$), we also can show that%
\begin{equation}%
\begin{cases}
||\partial_{x}^{l}L_{1}^{3}(n)||_{C_{0}([\sigma,T);L^{p_{3}}(\mathbb{R}^{2}%
))}\leq\epsilon,\ \ \forall\ \ 0\leq l\leq m-1,\\
||\partial_{x}^{m}L_{1}^{3}(n)||_{C_{\frac{1}{2}}((\sigma,T);L^{p_{3}%
}(\mathbb{R}^{2}))}\leq\epsilon.
\end{cases}
\label{Aux21.1}%
\end{equation}

Next, taking $K$ as in (\ref{AuxK}) and $(T-\sigma)$ small enough, we see that
$\Phi$ maps $B_{K}$ into itself. Moreover, we can show that $\Phi$ restricted
to $B_{K}$ is a contraction. Then, we have a fixed point $(n,c,\zeta)$ of
$\Phi$ in $B_{K}$, which is a solution of (\ref{CNS}). By the uniqueness
property in Theorem \ref{TLS}, we see that $(n,c,\zeta)$ coincides with the
original one in $(\sigma,T)$.

In particular, by the definition of $B_{K}$, we have proved that $(n,c,\zeta)$
satisfies
\[
\partial_{x}^{l}(n,c,\zeta)\in X_{1}(\sigma,T)\times X_{2}(\sigma,T)\times
X_{3}(\sigma,T),\text{ for }l=0,1,...,m.
\]
Take now $\sigma_{1}\in(\sigma,T)$. Set $T_{1}:=T-\sigma+\sigma_{1}$ and note
that (\ref{Aux21}) holds in the interval $(\sigma_{1},T_{1})$. By a similar
argument, we also can show that%

\[
\partial_{x}^{l}(n,c,\zeta)\in X_{1}(\sigma,T_{1})\times X_{2}(\sigma
,T_{1})\times X_{3}(\sigma,T_{1}),\text{ for }l=0,1,...,m,
\]
and, repeating this process a finite number of times, we arrive at%

\[
\partial_{x}^{l}(n,c,\zeta)\in X_{1}(\sigma,T_{0})\times X_{2}(\sigma
,T_{0})\times X_{3}(\sigma,T_{0}),\text{ for }l=0,1,...,m.
\]
Since $\sigma$ is arbitrary, it follows that
\begin{equation}
\partial_{x}^{l}(n,c,\zeta)\in\widetilde{X}_{1}(0,T_{0})\times\widetilde
{X}_{2}(0,T_{0})\times\widetilde{X}_{3}(0,T_{0}),\text{ for }l=0,1,...,m.
\label{KSreg00}%
\end{equation}

The case $m_{1}\geq1$ can be obtained from (\ref{KSreg00}). In fact, suppose
that (\ref{Aux8.1}) holds true for all $m_{1}\leq N-1$. We are going to show
that (\ref{Aux8.1}) is also valid for $m_{1}=N$. In fact, from system
(\ref{CNS}), note that the derivatives $\partial_{t}^{N}$ of $n,c,\zeta$ can
be expressed as sums of products of derivatives of $n,c,\zeta$ involving the
time partial derivative $\partial_{t}$ no more than $N-1$ times. But, by
induction hypothesis, those functions belong to the space $\widetilde{X}%
_{1}(0,T_{0})\times\widetilde{X}_{2}(0,T_{0})\times\widetilde{X}_{3}(0,T_{0}%
)$. Then, we can conclude that
\[
\partial_{t}^{m_{1}}\partial^{m_{2}}_{x}(n,c,\zeta)\in\widetilde{X}%
_{1}(0,T_{0})\times\widetilde{X}_{2}(0,T_{0})\times\widetilde{X}_{3}%
(0,T_{0}).
\]

Finally, using the previous case and proceeding similarly to the proof of
Proposition \ref{PropLp} and the above arguments, we obtain that
(\ref{Aux8.1}) is also valid for every $q_{1},q_{3}\in\lbrack1,\infty)$ and
$q_{2}\in\lbrack2,\infty),$ as desired.

\hfill\hfill$\diamond$\vskip12pt

From the proof of Proposition \ref{regularity}, note that we also have the
following result:

\begin{proposition}
\label{regularity-2}Let $(n,c,\zeta)$ be the mild solution of (\ref{CNS})
obtained in Theorem \ref{TLS}. Then, for $[t_{1},t_{2}]\subset(0,T_{0}),$
$q_{1},q_{3}\in\lbrack1,\infty),$ $q_{2}\in\lbrack2,\infty)$ and $m_{1}%
,m_{2}\in\mathbb{N}$ \ fixed, there exist constants $K_{n},$ $K_{c}$ and
$K_{\zeta}$, depending only on $q_{1},q_{2},q_{3}$,$t_{1},t_{2},m_{1},m_{2}$
and the norms of the initial data $n_{0},c_{0},\zeta_{0},$ such that for all
$n_{1}\leq m_{1}$, $n_{2}\leq m_{2}$ and $t\in\lbrack t_{1},t_{2}],$ we have
\begin{equation}
||\partial_{t}^{n_{1}}\partial_{x}^{n_{2}}n(t)||_{q_{1}}\leq K_{n}%
,\ ||\partial_{t}^{n_{1}}\partial_{x}^{n_{2}}c(t)||_{q_{2}}\leq K_{c},\text{
and }\ ||\partial_{t}^{n_{1}}\partial_{x}^{n_{2}}\zeta(t)||_{q_{3}}\leq
K_{\zeta}, \label{Aux8.2}%
\end{equation}
provided that $\nabla\phi\in H^{m_{2}}(\mathbb{R}^{2}).$
\end{proposition}

\section{Local-in-time solutions with measure data}

\label{sectionM} In the previous section, we obtain existence of solutions for
(\ref{CNS}) with initial data $\zeta_{0}$ and $n_{0}$ in $L^{1}(\mathbb{R}%
^{2})$, and $c_{0}\in L^{\infty}(\mathbb{R}^{2})$ with $\nabla c_{0}\in
L^{2}(\mathbb{R}^{2})$. In this part, inspired by \cite{Kato}, we are going to
obtain the existence of solution with the initial data $\zeta_{0}$ and $n_{0}$
being Radon measures, and preserving the same class for the data $c_{0}$. This
solution will be obtained as the limit of a sequence of solutions
$(n_{j},c_{j},\zeta_{j})$ with initial data $\zeta_{0,j}$ and $n_{0,j}$ in
$L^{1}(\mathbb{R}^{2})$.

For that, first recall that for each $\omega\in\mathcal{M}(\mathbb{R}^{2})$
there exists a sequence $\omega_{j}\in L^{1}(\mathbb{R}^{2})$ such that
$||\omega_{j}||_{1}\leq\pmb{|}\omega\pmb{|}$ and $\omega_{j}\rightarrow
\omega\in\mathcal{S^{\prime}}$. In other words, a Radon measure $\omega$ can
be approximated in $\mathcal{S^{\prime}}(\mathbb{R}^{2})$ by means of a
sequence of $L^{1}$-functions with $L^{1}$-norms not exceeding the total
variation of $\omega$. More precisely, it is sufficient to consider
$\varphi_{j}(x):=j^{2}\varphi(jx)$ and the sequence $\omega_{j}:=\varphi
_{j}\ast\omega$.

Next, consider $\zeta_{0},n_{0}\in\mathcal{M}(\mathbb{R}^{2})$ and $c_{0}\in
L^{\infty}(\mathbb{R}^{2})$ with $\nabla c_{0}\in L^{2}(\mathbb{R}^{2})$. By
the above comments, we have sequences $(\zeta_{0,j})_{j=1}^{\infty}$ and
$(n_{0,j})_{j=1}^{\infty}$ such that $\zeta_{0,j},n_{0,j}\in L^{1}%
(\mathbb{R}^{2})$,
\begin{equation}
||\zeta_{0,j}||_{1}\leq\pmb{|}\zeta_{0}\pmb{|}\text{ and }||n_{0,j}||_{1}%
\leq\pmb{|}n_{0}\pmb{|}, \label{aux-initial-bounded-1}%
\end{equation}
for all $j\in\mathbb{N}$, and $\zeta_{0,j}\rightarrow\zeta_{0}\in
\mathcal{S^{\prime}}$ and $n_{0,j}\rightarrow n_{0}\in\mathcal{S^{\prime}}$
when $j\rightarrow\infty$.

For each $j\in\mathbb{N}$, let $(n_{j},c_{j},\zeta_{j})$ be a solution of
(\ref{CNS}) in $(0,T_{j})$ with initial data $(n_{0,j},c_{0},\zeta_{0,j})$. By
Theorem \ref{TLS}, these solutions always exist for some appropriate $T_{j}$.
Using the feature of the approximation scheme, as well as the boundedness of
the initial data, we can obtain a $T_{0}>0$ such that $T_{0}\leq T_{j}$, for
all $j\in\mathbb{N}$.

The next lemma shows that we can obtain a subsequence of $(n_{j},c_{j}%
,\zeta_{j})_{j=1}^{\infty}$ that converge locally uniformly in $\mathbb{R}%
^{2}\times(0,T_{0})$ to a certain triple $(n,c,\zeta)$, where $T_{0}$ is the
same existence time in Theorem \ref{TLS}.

\begin{lemma}
\label{D7.2} Let $(n_{j},c_{j},\zeta_{j})_{j=1}^{\infty}$ defined as above,
$T_{0}>0$ as in Theorem \ref{TLS}, and $\nabla\phi\in H^{m+1}(\mathbb{R}^{2}%
)$, for $m\in\mathbb{N}$. Then, we can extract a subsequence (still denoted by
$(n_{j},c_{j},\zeta_{j})_{j=1}^{\infty}$) that converges locally uniformly in
$\mathbb{R}^{2}\times(0,T_{0})$ to a triple of $C^{m}$ functions $(n,c,\zeta)$.
\end{lemma}

\textit{Proof.} For each $i\in\mathbb{N}$, define the open set $K_{i}%
\subset\mathbb{R}^{2}\times(0,T_{0})$ by
\[
K_{i}:=(-i,i)^{2}\times\bigg(\frac{T_{0}}{i+1},\frac{T_{0}(i+1)}{i+2}\bigg).
\]
Notice that $K_{i}$ is bounded, $\underset{i\in\mathbb{N}}{\bigcup}%
K_{i}=\mathbb{R}^{2}\times(0,T_{0})$, and $K_{i}\subset K_{i+1}$. Also,
consider
\[
n_{j}^{i}:=n_{j}|_{K_{i}},\ \ c_{j}^{i}:=c_{j}|_{K_{i}},\text{ and}%
\ \ \zeta_{j}^{i}:=\zeta_{j}|_{K_{i}}.
\]

Let $l,k\in\mathbb{N}$, $q_{1}\in\lbrack1,\infty)$ and $q_{2}\in
\lbrack2,\infty).$ By Remark \ref{regularity} and (\ref{aux-initial-bounded-1}%
), using that $\nabla\phi\in H^{m+1}(\mathbb{R}^{2})$, the norms
\begin{equation}
||\partial_{t}^{k}\partial_{x}^{l}n_{j}(t)||_{q_{1}},||\partial_{t}%
^{k}\partial_{x}^{l}c_{j}(t)||_{q_{2}}\text{ and }||\partial_{t}^{k}%
\partial_{x}^{l}\zeta_{j}(t)||_{q_{1}} \label{aux-D7.2-1}%
\end{equation}
are bounded in any compact set $[t_{1},t_{2}]\subset(0,T_{0})$, uniformly
w.r.t. $j$, for all $0\leq\left\vert l\right\vert \leq m+1$. In particular,
taking $q_{1}=q_{2}=q\in(3,\infty)$, we have that $n_{j}^{i},c_{j}^{i}%
,\zeta_{j}^{i}\in W^{m+1,q}(K_{i})$. Moreover, $(\zeta_{j}^{i})_{j=1}^{\infty
}$, $(c_{j}^{i})_{j=1}^{\infty}$ and $(n_{j}^{i})_{j=1}^{\infty}$ are bounded
sequences in $W^{m+1,q}(K_{i})$. By compact Sobolev embedding, there exists
$\mathbb{N}_{1}\subset\mathbb{N}$ such that $(n_{j}^{1})_{j\in\mathbb{N}_{1}%
},$ $(c_{j}^{1})_{j\in\mathbb{N}_{1}}$ and $(\zeta_{j}^{1})_{j\in
\mathbb{N}_{1}}$ converge respectively to the functions $n^{1}:{K_{1}%
}\rightarrow\mathbb{R}$, $c^{1}:{K_{1}}\rightarrow\mathbb{R}$ and $\zeta
^{1}:{K_{1}}\rightarrow\mathbb{R}$ in $C^{m,\lambda}(\bar{K_{1}})$ with
$0<\lambda\leq1-3/q$.

Similarly, $(n_{j}^{2})_{j\in\mathbb{N}_{1}}$, $(c_{j}^{2})_{j\in
\mathbb{N}_{1}}$ and $(\zeta_{j}^{2})_{j\in\mathbb{N}_{1}}$ are bounded
sequences in $W^{m+1,q}(K_{2})$, and then there exists $\mathbb{N}_{2}%
\subset\mathbb{N}_{1}$ such that $(n_{j}^{2})_{j\in\mathbb{N}_{2}}$,
$(c_{j}^{2})_{j\in\mathbb{N}_{2}}$ and $(\zeta_{j}^{2})_{j\in\mathbb{N}_{2}}$
converge respectively to the functions $n^{2}:K_{2}\rightarrow\mathbb{R}$,
$c^{2}:K_{2}\rightarrow\mathbb{R}$ and $\zeta^{2}:K_{2}\rightarrow\mathbb{R}$
in $C^{m,\lambda}(\bar{K_{2}})$. By construction, we have that $n^{2}|_{K_{1}%
}\equiv n^{1}$, $c^{2}|_{K_{1}}\equiv c^{1}$ and $\zeta^{2}|_{K_{1}}%
\equiv\zeta^{1}$.

Proceeding in an inductive way, we can show that for each $i\in\mathbb{N}$,
there exists $\mathbb{N}_{i}\subset\mathbb{N}_{i-1}\subset...\subset
\mathbb{N}_{1}\subset\mathbb{N}$ such that the subsequence $(n_{j}^{i}%
)_{j\in\mathbb{N}_{i}},$ $(c_{j}^{i})_{j\in\mathbb{N}_{i}}$ and $(\zeta
_{j}^{i})_{j\in\mathbb{N}_{i}}$ converge to the functions $n^{i}%
:K_{i}\rightarrow\mathbb{R}$, $c^{i}:K_{i}\rightarrow\mathbb{R}$ and
$\zeta^{i}:K_{i}\rightarrow\mathbb{R}$ in $C^{m,\lambda}(\bar{K_{i}})$, respectively.

Next, define $\tilde{n}_{j}^{i}:\mathbb{R}^{2}\times(0,T_{0})\rightarrow
\mathbb{R}$ by
\[
\hspace{2cm}\tilde{n}_{j}^{i}(x):=%
\begin{cases}
n_{j}^{i}(x),\ \text{\ se\ }\ x\in K_{i},\\
0,\hspace{0.85cm}\ \text{\ otherwise,}%
\end{cases}
\]
$\tilde{c}_{j}^{i}:\mathbb{R}^{2}\times(0,T_{0})\rightarrow\mathbb{R}$ by
\[
\hspace{2cm}\tilde{c}_{j}^{i}(x):=%
\begin{cases}
c_{j}^{i}(x),\ \ \text{se}\ \ x\in K_{i},\\
0,\hspace{0.85cm}\ \ \text{otherwise,}%
\end{cases}
\]
and $\tilde{\zeta}_{j}^{i}:\mathbb{R}^{2}\times(0,T_{0})\rightarrow\mathbb{R}$
by
\[
\hspace{2cm}\tilde{\zeta}_{j}^{i}(x):=%
\begin{cases}
\zeta_{j}^{i}(x),\ \text{\ se}\ \ x\in K_{i},\\
0,\hspace{0.85cm}\text{\ \ otherwise.}%
\end{cases}
\]
Denoting $\mathbb{N}^{\ast}=(n_{1},n_{2},...),$ where $n_{i}$ is the $i$-th
term de $\mathbb{N}_{i}$, and employing a Cantor diagonal argument, we
conclude that the sequences $(\tilde{n}_{j}^{j})_{j\in\mathbb{N}^{\ast}}$,
$(\tilde{c}_{j}^{j})_{j\in\mathbb{N}^{\ast}}$ and $(\tilde{\zeta}_{j}%
^{j})_{j\in\mathbb{N}^{\ast}},$ as well as the corresponding $n_{j},c_{j}$ and
$\zeta_{j}$ (up to a subsequence), converge locally uniformly in
$\mathbb{R}^{2}\times(0,T_{0})$ to functions $n,c$ and $\zeta$ belonging to
$C^{m}\big(\mathbb{R}^{2}\times(0,T_{0})\big)$, respectively.\hfill
\hfill$\diamond$\vskip12pt

In order to conclude that the limit $(n,c,\zeta)_{j=1}^{\infty}$ of the
subsequence $(n_{j},c_{j},\zeta_{j})_{j=1}^{\infty}$ in Lemma \ref{D7.2} is
indeed a solution of (\ref{CNS}) for initial data $(n_{0},c_{0},\zeta_{0})$,
where $n_{0}$ and $\zeta_{0}$ are Radon measures, we need to show a suitable
equicontinuity property for $(n_{j},c_{j},\zeta_{j})_{j=1}^{\infty}$ and that
$(n,c,\zeta)_{j=1}^{\infty}$ satisfies the initial condition in the sense of distributions.

\begin{lemma}
\label{zetaequi} Let $(n_{j})_{j=1}^{\infty}$, $(c_{j})_{j=1}^{\infty}$ and
$(\zeta_{j})_{j=1}^{\infty}$ be the subsequences in Lemma \ref{D7.2} and let
$n,c,$ and $\zeta$ be their respective limits as $t\rightarrow0^{+}.$ Then,
$(n_{j})_{j=1}^{\infty}$, $(c_{j})_{j=1}^{\infty}$ and $(\zeta_{j}%
)_{j=1}^{\infty}$ are equicontinuous in $\mathcal{S^{\prime}}$ in the interval
$[0,T_{0})$. Moreover, $n(t)\rightarrow n_{0},$ $c(t)\rightarrow c_{0}$ and
$\zeta(t)\rightarrow\zeta_{0}$ in $\mathcal{S^{\prime}}$, as $t\rightarrow
0^{+}$.
\end{lemma}

\textit{Proof.} We shall show the statements only for the subsequence
$(n_{j})_{j=1}^{\infty}$ since the corresponding proofs for $(c_{j}%
)_{j=1}^{\infty}$ and $(\zeta_{j})_{j=1}^{\infty}$ are similar. Given $\psi
\in\mathcal{S}$ and denoting $u_{j}=S\ast\zeta_{j}$, we have that%

\begin{align*}
\partial_{t}\langle n_{j},\psi\rangle &  =\langle\Delta n_{j},\psi
\rangle-\langle\nabla\cdot(n_{j}u_{j}),\psi\rangle-\langle\nabla\cdot
(n_{j}\nabla c_{j}),\psi\rangle\\
&  =\langle n_{j},\Delta\psi\rangle+\langle n_{j}u_{j},\nabla\cdot\psi
\rangle+\langle n_{j}\nabla c_{j},\nabla\cdot\psi\rangle.
\end{align*}
Using \textbf{Condition A}, we can choose $q_{3},q,r,p_{1}^{\ast},q^{\ast
},r^{\ast}\in(1,\infty),$ such that $\frac{1}{q_{3}}=\frac{1}{p_{3}}-\frac
{1}{2}$, $\frac{1}{q}=\frac{1}{p_{1}}+\frac{1}{q_{3}}$, $\frac{1}{r}=\frac
{1}{p_{1}}+\frac{1}{p_{2}}$, $\frac{1}{p_{1}^{\ast}}=1-\frac{1}{p_{1}}$,
$\frac{1}{q^{\ast}}=1-\frac{1}{q}$ and $\frac{1}{r^{\ast}}=1-\frac{1}{r}$, and
then estimate
\begin{align*}
|\partial_{t}\langle n_{j},\psi\rangle|  &  \leq|\langle n_{j},\Delta
\psi\rangle|+|\langle n_{j}u_{j},\nabla\cdot\psi\rangle|+|\langle n_{j}\nabla
c_{j},\nabla\cdot\psi\rangle|\\
&  \leq||n_{j}||_{p_{1}}||\Delta\psi||_{p_{1}^{\ast}}+||n_{j}u_{j}%
||_{q}||\nabla\cdot\psi||_{q^{\ast}}+||n_{j}\nabla c_{j}||_{r}||\nabla
\cdot\psi||_{r^{\ast}}\\
&  \leq||n_{j}||_{p_{1}}||\Delta\psi||_{p_{1}^{\ast}}+||n_{j}||_{p_{1}}%
||u_{j}||_{q_{3}}||\nabla\cdot\psi||_{q^{\ast}}+||n_{j}||_{p_{1}}||\nabla
c_{j}||_{p_{2}}||\nabla\cdot\psi||_{r^{\ast}}.
\end{align*}
Applying (\ref{littlewood}) to the term $||u_{j}||_{q_{3}}$, we obtain
\[
|\partial_{t}\langle n_{j},\psi\rangle|\leq||n_{j}||_{p_{1}}||\Delta
\psi||_{p_{1}^{\ast}}+\sigma_{p_{3}}||n_{j}||_{p_{1}}||\zeta_{j}||_{p_{3}%
}||\nabla\cdot\psi||_{q^{\ast}}+||n_{j}||_{p_{1}}||\nabla c_{j}||_{p_{2}%
}||\nabla\cdot\psi||_{r^{\ast}}.
\]

Using Theorem \ref{TLS} together with the boundedness $||\zeta_{0,j}||_{1}%
\leq\pmb{|}\zeta_{0}\pmb{|}$ and $||n_{0,j}||_{1}\leq\pmb{|}n_{0}\pmb{|}$, we
arrive at $|\partial_{t}\langle n_{j},\psi\rangle|\leq C(t^{-\alpha_{1}%
}+t^{-\alpha_{1}-\alpha_{3}}+t^{-\alpha_{1}-\alpha_{2}})$, where $C>0$ is
independent of $t$. Since $\alpha_{1},(\alpha_{1}+\alpha_{2})$ and
$(\alpha_{1}+\alpha_{3})$ are less than $1$ (by (\ref{A4})), we have that
$C(t^{-\alpha_{1}}+t^{-\alpha_{1}-\alpha_{3}}+t^{-\alpha_{1}-\alpha_{2}})$ is
integrable on $[0,T_{0})$, and then we conclude the equicontinuity property.

In what follows, we deal with the convergence to the initial data. For
$\psi\in\mathcal{S}$ and $\epsilon>0$, we have that
\begin{equation}
|\langle n(t)-n_{0},\psi\rangle|\leq|\langle n(t)-n_{j}(t),\psi\rangle
|+|\langle n_{j}(t)-n_{0,j},\psi\rangle|+|\langle n_{0,j}-n_{0},\psi\rangle|.
\label{AuxEn}%
\end{equation}
Using the equicontinuity of $n_{j}$, there exists $j_{0}\in\mathbb{N}$ such
that
\begin{equation}
|\langle n(t)-n_{j}(t),\psi\rangle|<\frac{\epsilon}{3},\ \ \forall j\geq
j_{0}. \label{Auxn01}%
\end{equation}
Since $n_{0,j}\rightarrow n_{0}$ in $\mathcal{S^{\prime}}$, there exists
$j_{1}\in\mathbb{N}$ such that
\begin{equation}
|\langle n_{0,j}-n_{0},\psi\rangle|<\frac{\epsilon}{3},\ \ \forall j\geq
j_{1}. \label{Auxn02}%
\end{equation}
Moreover, since $n_{j}(t)\rightarrow n_{0,j}\in\mathcal{S^{\prime}}$ as
$t\rightarrow0$, there exists $T_{1}\in(0,T_{0})$ such that
\begin{equation}
|\langle n_{j}(t)-n_{0,j},\psi\rangle|<\frac{\epsilon}{3},\ \ \forall
t\in(0,T_{1}). \label{Auxn03}%
\end{equation}
Taking $j=\max\{j_{0},j_{1}\}$ in (\ref{AuxEn}) and using (\ref{Auxn01}),
(\ref{Auxn02}) and (\ref{Auxn03}), we obtain that
\[
|\langle n(t)-n_{0},\psi\rangle|<\epsilon,\ \ \forall t\in(0,T_{1}),
\]
which gives the desired convergence.\hfill\hfill$\diamond$\vskip12pt

Finally, collecting the properties obtained in Lemmas \ref{D7.2} and
\ref{zetaequi}, we have the following theorem:

\begin{theorem}
{(Measure initial data)} \label{TLSM} For $1\leq i\leq3$ and $m\in\mathbb{N}$,
let $p_{i},\alpha_{i}$ satisfy \textbf{Condition A }and let $\zeta_{0}%
,n_{0}\in\mathcal{M}(\mathbb{R}^{2})$ and $c_{0}\in L^{\infty}(\mathbb{R}%
^{2})$ with $\nabla c_{0}\in L^{2}(\mathbb{R}^{2})$. Assume also that
$\nabla\phi\in H^{2}(\mathbb{R}^{2}).$ If $||c_{0}||_{\infty}$ is sufficiently
small, then there exists a $T_{0}>0$ such that (\ref{CNS}) has a solution
$(n,c,\zeta)$ belonging to space $X$ (see (\ref{spacesolutions-0}%
)-(\ref{spacesolutions})), locally in $(0,T_{0}),$ with initial data
$(n_{0},c_{0},\zeta_{0})$, satisfying $(n,c,\zeta)\rightarrow(n_{0}%
,c_{0},\zeta_{0})$ in $\mathcal{S}^{\prime}$ as $t\rightarrow0^{+}.$
\end{theorem}

\begin{remark}
\label{rem-teo-measure}

\begin{itemize}
\item[$(i)$] (Further $L^{p}$-integrability and regularity) In view of Theorem
\ref{TLS} and Proposition \ref{PropLp}, the approximate sequence $(n_{j}%
,c_{j},\zeta_{j})_{j=1}^{\infty}$ in Lemma \ref{D7.2} belongs to space
(\ref{spacesolutions-0})-(\ref{spacesolutions}) for all $p_{1},p_{3}\in
\lbrack1,\infty)$, $p_{2}\in\lbrack2,\infty)$ and $\alpha_{1},\alpha
_{2},\alpha_{3}\geq0$ satisfying (\ref{A1}), where the corresponding norms in
(\ref{spacesolutions}) are bounded by initial data norms. Using
(\ref{aux-initial-bounded-1}), it follows that $(n_{j},c_{j},\zeta_{j}%
)_{j=1}^{\infty}$ is uniformly bounded in (\ref{spacesolutions-0}) and, by the
uniqueness of the limit in the sense of distributions, the solution
$(n,c,\zeta)$ also belongs to (\ref{spacesolutions-0}) under the above
conditions. Also, in view of Lemma \ref{D7.2} and its proof, we have that the
solution given in Theorem \ref{TLSM} is of class $C^{1}$ in $\mathbb{R}%
^{2}\times(0,T_{0})$, as well as satisfies (\ref{Aux8.1})-(\ref{Aux8.2}) with
$m_{1}\in\mathbb{N}$ and $m_{2}=2,$ since $\nabla\phi\in H^{2}(\mathbb{R}%
^{2}).$ In fact, we need this regularity of $\phi$ for performing the
approximation argument in the proof of Lemma \ref{D7.2} (see, e.g.,
(\ref{aux-D7.2-1})). Assuming that $\nabla\phi\in H^{m+1}(\mathbb{R}^{2})$
with $m\in\mathbb{N},$ we obtain that $(n,c,\zeta)$ is of class $C^{m}$ and
satisfies (\ref{Aux8.1})-(\ref{Aux8.2}) for all $m_{1}\in\mathbb{N}$ and
$m_{2}=m+1.$

\item[$(ii)$] (Uniqueness of solutions) Our existence result considers initial
measure $n_{0}$ and $\zeta_{0}$ without assuming smallness conditions on them.
However, since we construct the solutions via a limit of approximate
solutions, we have lost the uniqueness property. For $\zeta_{0},n_{0}\in
L^{1}(\mathbb{R}^{2}),$ the uniqueness follows from Theorem \ref{TLS}. For a
general $f\in\mathcal{M}(\mathbb{R}^{2}),$ we have the decomposition
$f=f_{a}+f_{c}$ in its atomic and continuous parts and
$\pmb{|}f\pmb{|}=\pmb{|}f_{a}\pmb{|}+\pmb{|}f_{c}\pmb{|}$. Recalling the
seminorm (\ref{seminorm1}), by \cite[p.951]{Kato},\cite[p.245]{Giga1}, it
follows that
\begin{equation}
\pmb{\lVert}e^{t\Delta}f\pmb{\lVert}_{p,\alpha}\leq\sigma_{p}\pmb{|}f_{a}%
\pmb{|}, \label{est-sing-1}%
\end{equation}
where $\pmb{|}f_{a}\pmb{|}$ is the total variation of the atomic part of $f$.
So, $\pmb{\lVert}e^{t\Delta}f\pmb{\lVert}_{p,\alpha}=0$ when $f$ is a
continuous measure (i.e., $f_{a}=0$). Denoting $d=(d_{1},d_{2},...)\in
l^{1}(\mathbb{N})$ and writing $f_{a}=\sum_{j=1}^{\infty}d_{j}\delta
(x-a_{j}),$ we have that
\begin{equation}
\pmb{|}f_{a}\pmb{|}=\sum_{j=1}\left\vert d_{j}\right\vert <\infty.
\label{est-sing-2}%
\end{equation}
For $\zeta_{0},n_{0}\in\mathcal{M}(\mathbb{R}^{2}),$ decomposing $\zeta
_{0}=\zeta_{0,a}+\zeta_{0,c}$ and $n_{0}=n_{0,a}+n_{0,c}$ and using
(\ref{est-sing-1}), we can make $|||e^{t\Delta}n_{0}|||_{p_{1},\alpha_{1}}$
and $|||e^{t\Delta}\zeta_{0}|||_{p_{3},\alpha_{3}}$ small in
(\ref{aux-data-epsilon-1}) provided that $\pmb{|}n_{0,a}\pmb{|}$ and
$\pmb{|}\zeta_{0,a}\pmb{|}$ are small enough, and then recover the
existence-uniqueness statement in Theorem \ref{TLS}. In particular, those
conclusions hold true for continuous measures, that is, without atomic part
(i.e., $n_{0,a}=\zeta_{0,a}=0$). Note that filaments in $\mathbb{R}^{2}$
(i.e., measures concentrated on smooth curves) are continuous measures and
then we have the uniqueness property for this kind of initial data. Finally,
in view of (\ref{est-sing-2}), we can express $\pmb{|}n_{0,a}\pmb{|}=\sum
_{j=1}^{\infty}\left\vert d_{n_{j}}\right\vert $ and $\pmb{|}\zeta
_{0,a}\pmb{|}=\sum_{j=1}^{\infty}\left\vert d_{\zeta_{j}}\right\vert $ for
suitable coefficients $(d_{n_{j}}),(d_{\zeta_{j}})\in l^{1}(\mathbb{N}),$
which allow to provide uniqueness-conditions (maybe explicitly) in terms of
the breadth of values $d_{n_{j}}$ and $d_{\zeta_{j}}$ of the Diracs
$\delta(x-a_{n_{j}})$ and $\delta(x-a_{\zeta_{j}})$, respectively. In this
direction, we remark that Bedrossian and Masmoudi \cite{Bedro-Masmoudi} showed
local-in-time existence-uniqueness for the Patlak-Keller-Segel with initial
data $\mu\in\mathcal{M}(\mathbb{R}^{2})$ satisfying $\max_{x\in\mathbb{R}^{2}%
}\mu(x)<8\pi$.

\item[$(iii)$] (Persistence property) In view of Remark \ref{Rem-RTL-PropLp}
$(i)$, Lemma \ref{zetaequi} and Theorem \ref{TLSM}, considering $p_{1}%
=p_{3}=1$ and $\alpha_{1}=\alpha_{3}=0$, we obtain a solution $(n,c,\zeta)$
such that $n\in BC((0,T_{0});\mathcal{M}(\mathbb{R}^{2}))$, $c\in
BC((0,T_{0});L^{\infty}(\mathbb{R}^{2}))$ and $\zeta\in BC((0,T_{0}%
);\mathcal{M}(\mathbb{R}^{2}))$ with $n(t)\rightarrow n_{0},$ $c(t)\rightarrow
c_{0}$ and $\zeta(t)\rightarrow\zeta_{0}$ in $\mathcal{S^{\prime}}%
(\mathbb{R}^{2})$ as $t\rightarrow0^{+}$, that is, the time-continuity at
$t=0^{+}$ is taken in the sense of tempered distributions. Also, note that the
initial measures $n_{0},\zeta_{0}$ can be arbitrarily large. Moreover, for
each $t\in(0,T_{0}),$ we have a regularizing effect and $n(\cdot
,t),\zeta(\cdot,t)\in L^{1}(\mathbb{R}^{2}),$ that is, they are absolutely
continuous measures in positive times.


\end{itemize}
\end{remark}

\section{Global-in-time solutions}

\label{sectionG} In this section we show that the solutions constructed in
Theorem \ref{TLS} and Theorem \ref{TLSM} are global in time. For that, we are
going to show some $L^{p}$-boundedness uniformly in time which are the subject
of the next three propositions. In them, we should implicitly assume at least
the conditions of Theorem \ref{TLS}, Propositions \ref{PropLp} and
\ref{regularity} as well as employ the properties contained in those results.

Before proceeding, by the first equation in (\ref{CNS}), we know that $n\geq0$
preserves its mass, i.e.,
\begin{equation}
||n(t)||_{1}=||n(\sigma)||_{1},\text{ for all }t\in\lbrack\sigma,T_{max}),
\label{aux-mass-n-1}%
\end{equation}
where $\sigma\geq0$. Moreover, by maximum principle, we have that
\begin{equation}
\left\Vert c(t)\right\Vert _{\infty}\leq\left\Vert c(\sigma)\right\Vert
_{\infty},\text{ for all }t\in\lbrack\sigma,T_{\max}). \label{aux-Linfty-c-1}%
\end{equation}

The next proposition provides an estimate of the $L^{p}$-norm of $n$ under a
smallness condition on $||c(\sigma)||_{\infty}$. The proof can be obtained by
following closely the arguments in \cite[Proposition 2]{Chae2} and then we
omit it. However, for the convenience of the reader, let us observe some
necessary care in the adaptation process. In \cite[Proposition 2]{Chae2} the
authors considered solutions (see also \cite{Chae})%
\begin{equation}
(n,c,u)\in L^{\infty}((0,T);H^{m-1}(\mathbb{R}^{2})\times H^{m}(\mathbb{R}%
^{2})\times H^{m}(\mathbb{R}^{2})),\text{ }m\geq3, \label{aux-space-1010-2}%
\end{equation}
where $T>0$ is an existence time and $(n(\cdot,t),c(\cdot,t),u(\cdot,t))$ is
also continuous in time $t$. In particular, their solutions satisfy
\begin{equation}
(n(\cdot,t),c(\cdot,t),u(\cdot,t))\in H^{m-1}(\mathbb{R}^{2})\times
H^{m}(\mathbb{R}^{2})\times H^{m}(\mathbb{R}^{2})\text{ with }m\geq3,
\label{aux-space-1010}%
\end{equation}
for each $t>0$. Note that in (\ref{aux-space-1010}), the velocity
$u(\cdot,t)\in H^{m}(\mathbb{R}^{2}),$ and then $u(\cdot,t)\in L^{2}%
(\mathbb{R}^{2}),$ for each $t>0$, which is not verified by our solution. In
fact, using that $\zeta(\cdot,t)\in L^{q_{3}}(\mathbb{R}^{2})$ with $q_{3}$ as
in Theorem \ref{TLS} and Proposition \ref{PropLp} and applying the Biot-Savart
law (\ref{velocity1}) and HLS inequality (\ref{littlewood}), we get only that
$u(\cdot,t)\in L^{q}(\mathbb{R}^{2})$ with $q>2$, for each $t>0$, likewise for
its derivatives. One of the main reasons that led to this is that we are
considering initial vorticity $\zeta_{0}$ in $L^{1}(\mathbb{R}^{2})$ or, more
generally, $\zeta_{0}$ in $\mathcal{M}(\mathbb{R}^{2})$. Moreover, they
assumed $\nabla^{l}\phi\in L^{\infty}(\mathbb{R}^{2})$ for all $1\leq
\left\vert l\right\vert \leq m,$ while we are going to consider $\nabla\phi\in
H^{2}(\mathbb{R}^{2}).$ Therefore, the arguments in \cite[Proposition
2]{Chae2} need to be adapted using the level of $L^{p}$-integrability verified
by our solution and its derivatives according to Section \ref{sectionL}.

\begin{proposition}
\label{Nboundedness} Let $p\in(1,\infty)$ and let $n$ be as in Theorem
\ref{TLS} with the improvement in Proposition \ref{PropLp} and the maximal
existence time $T_{\max}$. Assume also that $n(\sigma)\in L^{p}(\mathbb{R}%
^{2})$ for some $\sigma\in(0,T_{max}).$ If $||c(\sigma)||_{\infty}\leq\frac
{1}{24p}$, then
\begin{equation}
||n(t)||_{p}\leq C_{n,p},\ \ \ \ \text{for all }t\in\lbrack\sigma,T_{\max}),
\label{C_p}%
\end{equation}
where the constant $C_{n,p}>0$ depends only on $p,\left\Vert n(\sigma
)\right\Vert _{p}$ and $||c(\sigma)||_{\infty}$.
\end{proposition}

In the sequel we obtain a suitable control for the $L^{p}$-norm of $\zeta$,
uniformly in time.

\begin{proposition}
\label{zetaboundedness} Let $p\in(1,\infty)$ and let $\zeta$ be as in Theorem
\ref{TLS} with the improvement in Proposition \ref{PropLp} and the maximal
existence time $T_{\max}$. Suppose also that $n(\sigma),\zeta(\sigma)\in
L^{4}(\mathbb{R}^{2})\cap L^{p}(\mathbb{R}^{2})$, $\nabla\phi\in L^{\infty
}(\mathbb{R}^{2})$ and $||c(\sigma)||_{\infty}\leq\min\{\frac{1}{24p},\frac
{1}{96}\}$, for some $\sigma\in(0,T),$ where $T\leq T_{\max}$. If $T<\infty$,
then
\begin{equation}
||\zeta(t)||_{p}\leq C_{\zeta,p},\ \ \ \text{for all }t\in\lbrack\sigma,T),
\label{Aux09}%
\end{equation}
where the constant $C_{\zeta,p}>0$ depends on $p,C_{n,p},C_{n,4},\left\Vert
\nabla\phi\right\Vert _{L^{\infty}},||\zeta(\sigma)||_{1},\left\Vert
\zeta(\sigma)\right\Vert _{L^{4}\cap L^{p}}$ and $||c(\sigma)||_{\infty}$, as
well as continuously on $T.$
\end{proposition}

\textit{Proof.} We first show (\ref{Aux09}) for $p\in(2,\infty)$. Indeed,
testing the third equation in (\ref{CNS}) with $\zeta|\zeta|^{p-2}$, we
obtain
\begin{equation}
\frac{1}{p}\frac{d}{dt}\int|\zeta|^{p}+\int(u\cdot\nabla\zeta)\zeta
|\zeta|^{p-2}-\int(\Delta\zeta)\zeta|\zeta|^{p-2}=-\int\nabla^{\perp}%
\cdot(n\nabla\phi)\zeta|\zeta|^{p-2}. \label{aux-vort-1}%
\end{equation}
Integrating by parts and using that $\nabla\cdot u=0$, we arrive at
\[
\int(u\cdot\nabla\zeta)\zeta|\zeta|^{p-2}=\frac{1}{p}\int u\cdot\nabla
(|\zeta|^{p})=-\frac{1}{p}\int(\nabla\cdot u)|\zeta|^{p}=0.
\]
Moreover, we have that
\[
-\int(\Delta\zeta)\zeta|\zeta|^{p-2}=\int(\nabla\zeta)\cdot\nabla(\zeta
|\zeta|^{p-2})=(p-1)\int|\nabla\zeta|^{2}|\zeta|^{p-2}%
\]
and
\begin{equation}
-\int\nabla^{\perp}\cdot(n\nabla\phi)\zeta|\zeta|^{p-2}=\int(n\nabla
\phi)\nabla^{\perp}(\zeta|\zeta|^{p-2})=(p-1)\int(n\nabla\phi)(\nabla^{\perp
}\zeta)|\zeta|^{p-2}. \label{Aux08}%
\end{equation}
Using Cauchy-Schwartz inequality in (\ref{Aux08}), we obtain that
\[
-\int\nabla^{\perp}\cdot(n\nabla\phi)\zeta|\zeta|^{p-2}\leq\frac{(p-1)}{4}%
\int(n\nabla\phi)^{2}|\zeta|^{p-2}+(p-1)\int|\nabla^{\perp}\zeta|^{2}%
|\zeta|^{p-2}.
\]
Noting that $|\nabla^{\perp}\zeta|^{2}=|\nabla\zeta|^{2}$, it follows that
\[
\frac{1}{p}\frac{d}{dt}\int|\zeta|^{p}\leq\frac{(p-1)}{4}\int(n\nabla\phi
)^{2}|\zeta|^{p-2}.
\]
If $p>2$, define $q:=\frac{p}{p-2}$. Note that if $q_{1}:=\frac{p}{2}$, then
$\frac{1}{q}+\frac{1}{q_{1}}=1$. Young inequality yields
\[
(n\nabla\phi)^{2}|\zeta|^{p-2}\leq\frac{2}{p}(n\nabla\phi)^{p}+\frac{p-2}%
{p}|\zeta|^{p}%
\]
and then
\begin{align*}
\frac{d}{dt}\int|\zeta|^{p}  &  \leq\frac{2(p-1)}{4}\int(n\nabla\phi
)^{p}+\frac{(p-1)(p-2)}{4}\int|\zeta|^{p}\\
&  \leq\frac{2(p-1)}{4}||\nabla\phi||_{\infty}^{p}C_{n,p}^{p}+\frac
{(p-1)(p-2)}{4}\int|\zeta|^{p}\ \ ,
\end{align*}
for all $t\in\lbrack\sigma,T)$, where above we have used (\ref{C_p}). It
follows from Gronwall inequality that
\[
\int|\zeta|^{p}\leq e^{t\big(\frac{(p-1)(p-2)}{4}\big)}\bigg(e^{-\sigma
\big(\frac{(p-1)(p-2)}{4}\big)}||\zeta(\sigma)||_{p}^{p}+(t-\sigma
)\bigg(\frac{2(p-1)}{4}||\nabla\phi||_{\infty}^{p}C_{n,p}^{p}\bigg)\bigg).
\]
Thus, for $T<\infty$,
\begin{equation}
||\zeta(t)||_{p}\leq C_{\zeta,p},\text{ for all }t\in\lbrack\sigma,T),
\label{Aux10}%
\end{equation}
where
\[
C_{\zeta,p}^{p}:=e^{T\big(\frac{(p-1)(p-2)}{4}\big)}\bigg(||\zeta
(\sigma)||_{p}^{p}+T\bigg(\frac{2(p-1)}{4}||\nabla\phi||_{\infty}^{p}%
C_{n,p}^{p}\bigg)\bigg).
\]

For $p\in(1,2]$, we can estimate the $L^{p}$-norm of $\zeta(t)$ by
interpolation. Indeed, since $||\zeta(t)||_{1}\leq||\zeta(\sigma)||_{1},$ for
all $t\in\lbrack\sigma,T)$, taking $\theta_{p}=\frac{4-p}{3p}$ and
$C_{\zeta,4}$ (case $p=4$) as in (\ref{Aux10}), we have
\begin{equation}
||\zeta(t)||_{p}\leq||\zeta(t)||_{1}^{\theta_{p}}||\zeta(t)||_{4}%
^{(1-\theta_{p})}\leq||\zeta(\sigma)||_{1}^{\theta_{p}}C_{\zeta,4}%
^{(1-\theta_{p})}, \label{Aux-11-11}%
\end{equation}
for all $t\in\lbrack\sigma,T).$ The result follows from (\ref{Aux10}) and
(\ref{Aux-11-11}).\hfill\hfill$\diamond$\vskip12pt

Until the moment, we have a uniform control of the $L^{p}$-norm of $n$ and
$\zeta$. Also, we will need a uniform control of the $L^{p}$-norm of $\nabla
c$. This is the subject of the proposition below.

\begin{proposition}
\label{GCboundedness} Let $p\in(2,\infty),$ $q\in(1,2),$ and let $c$ as in
Theorem \ref{TLS} with the improvement in Proposition \ref{PropLp} and the
maximal existence time $T_{\max}$. Suppose also that $n(\sigma)\in
L^{2}(\mathbb{R}^{2}),$ $\zeta(\sigma)\in L^{q}(\mathbb{R}^{2}),$ $\nabla
c(\sigma)\in L^{p}(\mathbb{R}^{2}),$ and $||c(\sigma)||_{\infty}\leq\frac
{1}{96}$, for some $\sigma\in(0,T),$ where $T\leq T_{\max}$. If $T<\infty$,
then
\begin{equation}
||\nabla c(t)||_{p}\leq C_{c,p},\ \ \ \ \forall t\in\lbrack\sigma,T),
\label{GC-est}%
\end{equation}
where the constant $C_{c,p}>0$ depends continuously on $T,$ as well as on
$p,q,C_{n,2},C_{\zeta,q},\left\Vert \nabla c(\sigma)\right\Vert _{p}$ and
$||c(\sigma)||_{\infty}$. The constants $C_{n,2}$ and $C_{\zeta,q}$ are as in
Propositions \ref{Nboundedness} and \ref{zetaboundedness}, respectively.
\end{proposition}

\textit{Proof.} By the integral formulation (\ref{IF}) of $c$, but starting
from $\sigma$ instead of $0$, we see that
\begin{equation}
\nabla c(t)=e^{(t-\sigma)\Delta}\nabla c(\sigma)-\int_{\sigma}^{t}\nabla
e^{(t-s)\Delta}(nc)(s)ds-\int_{\sigma}^{t}\nabla e^{(t-s)\Delta}((S\ast
\zeta)\cdot\nabla c)(s)ds. \label{gradc-aux-1}%
\end{equation}
By heat semigroup properties, it follows that
\begin{equation}
\left\Vert e^{(t-\sigma)\Delta}(\nabla c(\sigma))\right\Vert _{p}%
\leq\left\Vert \nabla c(\sigma)\right\Vert _{p}. \label{Aux05}%
\end{equation}
For the second term in the R.H.S. of (\ref{gradc-aux-1}), using
(\ref{est-sg-1}), (\ref{C_p}), and $\left\Vert c(t)\right\Vert _{\infty}%
\leq\left\Vert c(\sigma)\right\Vert _{\infty}$, we can estimate%

\begin{align}
\bigg|\bigg|\int_{\sigma}^{t}\nabla e^{(t-s)\Delta}(nc)(s)ds\bigg|\bigg|_{p}
&  \leq C\int_{\sigma}^{t}(t-s)^{-\frac{1}{2}-({\frac{1}{2}}-{\frac{1}{p}}%
)}||nc(s)||_{2}ds\nonumber\\
&  \leq C\int_{\sigma}^{t}(t-s)^{-\frac{1}{2}-(\frac{1}{2}-\frac{1}{p}%
)}||n(s)||_{2}||c(s)||_{\infty}ds\nonumber\\
&  \leq C_{n,2}||c(\sigma)||_{\infty}\int_{\sigma}^{t}(t-s)^{-\frac{1}%
{2}-({\frac{1}{2}}-{\frac{1}{p}})}ds\hspace{2cm}\nonumber\\
&  \leq CT^{1/p}\left\Vert c(\sigma)\right\Vert _{\infty}. \label{Aux05-2}%
\end{align}
For the third term, consider $q_{1}\in\lbrack1,p)$ and $r\in(2,\infty)$ such
that $\frac{1}{q_{1}}=\frac{1}{p}+\frac{1}{r}$ and $\frac{1}{r}=\frac{1}%
{q}-\frac{1}{2}$. Then, using estimates (\ref{est-sg-1}), (\ref{littlewood})
and (\ref{Aux09}), we obtain that%

\begin{align}
\bigg|\bigg|\int_{\sigma}^{t}\nabla e^{(t-s)\Delta}((S\ast\zeta)\cdot\nabla
c)(s)ds\bigg|\bigg|_{p}  &  \leq C\int_{\sigma}^{t}(t-s)^{-\frac{1}{2}%
-(\frac{1}{q_{1}}-\frac{1}{p})}||(S\ast\zeta)(s)||_{r}||\nabla c(s)||_{p}%
ds\ \ \ \ \nonumber\\
&  \leq C\sigma_{q}\int_{\sigma}^{t}(t-s)^{-\frac{1}{2}-\frac{1}{r}}%
||\zeta(s)||_{q}||\nabla c(s)||_{p}ds\ \ \nonumber\\
&  \leq C\sigma_{q}C_{\zeta,q}\int_{\sigma}^{t}(t-s)^{-\frac{1}{2}-\frac{1}%
{r}}||\nabla c(s)||_{p}ds.\ \label{Aux05-3}%
\end{align}

Taking $L_{\sigma}(T)=\left\Vert \nabla c(\sigma)\right\Vert _{p}%
+CT^{1/p}\left\Vert c(\sigma)\right\Vert _{\infty}$ and putting together
(\ref{Aux05}), (\ref{Aux05-2}) and (\ref{Aux05-3}) yield%
\[
\left\Vert \nabla c(t)\right\Vert _{p}\leq L_{\sigma}(T)+C\int_{\sigma}%
^{t}(t-s)^{-\delta}||\nabla c(s)||_{p}ds,
\]
for all $t\in\lbrack\sigma,T),$ with $\delta=\frac{1}{2}+\frac{1}{r}<1$. Now,
using fractional Gronwall-type inequalities (see, e.g., \cite{Webb-1}), we get
estimate (\ref{GC-est}).\hfill\hfill$\diamond$\vskip12pt

With the above boundedness in hand, we are in position to prove that the
solutions obtained in Theorem \ref{TLS} are global in time.

\begin{theorem}
\label{global}Assume the hypotheses of Theorem \ref{TLS} with the improvement
in Proposition \ref{PropLp} and let $(n,c,\zeta)$ be the unique solution for
(\ref{CNS}) obtained in that theorem with initial data $(n_{0},c_{0},\zeta
_{0})$. Suppose also that $\nabla\phi\in H^{2}(\mathbb{R}^{2}),$
$||c_{0}||_{\infty}\leq\delta_{0}$ for some suitable small $\delta_{0},$ and
let $T_{max}$ be its maximal existence time. Then, $T_{max}=\infty$. Moreover,
for each $T\in(0,\infty),$ we have that%
\begin{equation}
(n,c,\zeta)\in X=X_{1}\times X_{2}\times X_{3}, \label{space-measure-aux-1}%
\end{equation}
where
\begin{equation}%
\begin{split}
X_{1}=  &  C_{\alpha_{1}}((0,T);L^{p_{1}}(\mathbb{R}^{2})),\\
X_{2}=  &  \{c\in C_{0}((0,T);L^{\infty}(\mathbb{R}^{2}):\nabla c\in
{C_{\alpha_{2}}((0,T);L^{p_{2}}(\mathbb{R}^{2}))}\},\\
X_{3}=  &  C_{\alpha_{3}}((0,T);L^{p_{3}}(\mathbb{R}^{2})),
\end{split}
\label{space-measure-aux-2}%
\end{equation}
with $p_{1},p_{3}\in\lbrack1,\infty)$, $p_{2}\in\lbrack2,\infty)$ and
$\alpha_{1},\alpha_{2},\alpha_{3}\geq0$ satisfying (\ref{A1}).
\end{theorem}

\textit{Proof.} Let $(n,c,\zeta)$ be the solution of Theorem \ref{TLS} with
initial data $(n_{0},c_{0},\zeta_{0}).$ For positive times, since $\nabla
\phi\in H^{2}(\mathbb{R}^{2})\subset L^{\infty}(\mathbb{R}^{2}),$ this
solution has the needed regularity and properties in order to employ
Propositions \ref{Nboundedness}, \ref{zetaboundedness} and \ref{GCboundedness}
(see, e.g., Propositions \ref{PropLp}, \ref{regularity} and \ref{regularity-2}).

Suppose by contradiction that $T_{max}<\infty$. For any $\sigma\in(0,T_{max}%
)$, consider $\tilde{n}_{0}:=n(\sigma)$, $\tilde{c}_{0}:=c(\sigma)$ and
$\tilde{\zeta}_{0}:=\zeta(\sigma),$ and the corresponding solution
$(n_{\sigma},c_{\sigma},\zeta_{\sigma})$ of (\ref{CNS}) in $(0,T_{\sigma})$
obtained from Theorem \ref{TLS} with initial data $(\tilde{n}_{0},\tilde
{c}_{0},\tilde{\zeta}_{0}).$

Next, for a moment, take fixed indexes $\tilde{p}_{1},\tilde{p}_{3}\in
\lbrack4/3,2)$, $\tilde{p}_{2}\in\lbrack2,4)$ and $\tilde{\alpha}_{1}%
,\tilde{\alpha}_{2},\tilde{\alpha}_{3}>0$ satisfying (\ref{A1}). Then,
recalling the properties of the solution $(n,c,\zeta)$ (see Propositions
\ref{PropLp} and \ref{regularity}) and using Propositions \ref{Nboundedness},
\ref{zetaboundedness} and \ref{GCboundedness}, for $\tilde{T}\in
(\sigma,T_{\max}]$ we have the estimate
\begin{equation}
||n(\cdot,t)||_{\tilde{p}_{1}}+||c(\cdot,t)||_{\infty}+||\nabla c(\cdot
,t)||_{\tilde{p}_{2}}+||\zeta(\cdot,t)||_{\tilde{p}_{3}}\leq C,\ \ \forall
t\in(\sigma,\tilde{T}), \label{BU}%
\end{equation}
where $C>0$ depends continuously on $\tilde{T}$ and a finite set of $L^{p}%
$-norms of the local-in-time solution $(n,c,\zeta)$ (and $\nabla c$) at fixed
time $\sigma.$ Here, note that we need to take $\delta_{0}=\min\{\frac
{1}{24\tilde{p}_{1}},\frac{1}{96}\}=\frac{1}{96}$ according to the
above-mentioned propositions.

Consider $\tilde{T}=T_{\max}<\infty$ and note that, by boundedness (\ref{BU})
and Remark \ref{Rem-RTL-PropLp} $(ii)$, we can obtain an existence time $T>0$
with $T<T_{\tilde{\sigma}}$, for all $\tilde{\sigma}\in\lbrack\sigma,T_{\max
})$. Set $\sigma_{0}:=\max\{\sigma,T_{max}-\frac{T}{2}\}$. By Theorem
\ref{TLS} and Remark \ref{Rem-RTL-PropLp} $(ii)$, there exists a solution
$(\tilde{n},\tilde{c},\tilde{\zeta})$ of (\ref{CNS}) in $(0,T)$ starting from
$(n({\sigma_{0}}),c({\sigma_{0}}),\zeta({\sigma_{0}}))$. But, by the
uniqueness statement in Theorem \ref{TLS}, we know that $(n,c,\zeta)$
coincides with $(\tilde{n},\tilde{c},\tilde{\zeta})$ in $(\sigma_{0},T_{max}%
)$. More precisely, for all $t\in(0,T)$, we have that
\[
\tilde{n}(t)=n(t+\sigma_{0}),\ \ \tilde{c}(t)=c(t+\sigma_{0}),\text{
and}\ \ \tilde{\zeta}(t)=\zeta(t+\sigma_{0}).
\]

So, defining
\[
{\hat{n}}(t):=%
\begin{cases}
n(t),\ \ \text{if}\ \ t\in(0,\sigma_{0}),\\
\tilde{n}(t-\sigma_{0}),\ \ \text{if}\ \ t\in\lbrack\sigma_{0},T_{max}+T/2),
\end{cases}
\]%
\[
{\hat{c}}(t):=%
\begin{cases}
c(t),\ \ \text{if}\ t\in(0,\sigma_{0}),\\
\tilde{c}(t-\sigma_{0}),\ \ \text{if}\ \ t\in\lbrack\sigma_{0},T_{max}+T/2),
\end{cases}
\]
and%
\[
{\hat{\zeta}}(t):=%
\begin{cases}
\zeta(t),\ \ \text{if}\ \ t\in(0,\sigma_{0}),\\
\tilde{\zeta}(t-\sigma_{0}),\ \ \text{if}\ \ t\in\lbrack\sigma_{0}%
,T_{max}+T/2),
\end{cases}
\]
we see that $(\hat{n},\hat{c},\hat{\zeta})$ is a solution of (\ref{CNS}) in
$(0,T_{max}+\frac{T}{2})$ with initial data $(n_{0},c_{0},\zeta_{0}),$ and
$(\hat{n},\hat{c},\hat{\zeta})$ coincides with $(n,c,\zeta)$ in $(0,T_{max})$.
But this contradicts the definition of $T_{max}$ and then $T_{max}=\infty$. It
follows that the solution belongs to (\ref{space-measure-aux-1}) with the
fixed $\tilde{p}_{1},\tilde{p}_{2},\tilde{p}_{3},\tilde{\alpha}_{1}%
,\tilde{\alpha}_{2},\tilde{\alpha}_{3}$, for all finite $T>0$. Furthermore,
using the condition on $\tilde{p}_{1},\tilde{p}_{2},\tilde{p}_{3}$ and
performing a bootstrapping argument, similarly to Proposition \ref{PropLp}, we
obtain that $(n,c,\zeta)$ belongs to (\ref{space-measure-aux-1}) for every
$p_{1},p_{3}\in\lbrack1,\infty)$, $p_{2}\in\lbrack2,\infty)$ and $\alpha
_{1},\alpha_{2},\alpha_{3}\geq0$ satisfying (\ref{A1}).\hfill\hfill$\diamond$\vskip12pt

\begin{remark}
\label{Rem-aux-borderline}For example, in the proof above, the borderline case
$p_{1}=p_{3}=1,$ $p_{2}=2$ and $\alpha_{1}=\alpha_{2}=\alpha_{3}=0$ can be
reached by estimating the operators $B_{i,j}^{k}$ and $L_{1}^{3}$ via the norm
of (\ref{space-measure-aux-1}) with $p_{1}=p_{3}=4/3$ and $p_{2}=8/3.$
\end{remark}

\begin{theorem}
\label{global-M}Under the hypotheses of Theorem \ref{TLSM} with the
improvement in Remark \ref{rem-teo-measure} $(i)$. Let $(n,c,\zeta)$ be the
solution for (\ref{CNS}) obtained in that theorem with $\nabla\phi\in
H^{2}(\mathbb{R}^{2})$, initial data $\zeta_{0},n_{0}\in\mathcal{M}%
(\mathbb{R}^{2})$ and $c_{0}\in L^{\infty}(\mathbb{R}^{2})$ with $\nabla
c_{0}\in L^{2}(\mathbb{R}^{2}).$ Assume additionally that $||c_{0}||_{\infty
}\leq\delta_{0}$ for some suitable small $\delta_{0}$. Then, $(n,c,\zeta)$ can
be extended globally in time and belongs to the class
(\ref{space-measure-aux-1})-(\ref{space-measure-aux-2}), for every
$T\in(0,\infty).$ The solution is unique provided that the atomic parts of the
measures $n_{0}$ and $\zeta_{0}$ are small enough, such as in the case of
atomless measures.
\end{theorem}

\textit{Proof.} First note that $\nabla\phi\in H^{2}(\mathbb{R}^{2})\subset
L^{\infty}(\mathbb{R}^{2}).$ Let $(n,c,\zeta)$ be the solution constructed in
Theorem \ref{TLSM}. Recalling also their properties in Remark
\ref{rem-teo-measure} $(i)$ and $(iii)$, note that we can employ Propositions
\ref{Nboundedness}, \ref{zetaboundedness} and \ref{GCboundedness}. Therefore,
considering a suitable $\sigma>0$ and a start point $(n(\sigma),c(\sigma
),\zeta(\sigma))$, we can proceed as in the extension argument in the proof of
Theorem \ref{global} provided that $\delta_{0}=\min\{\frac{1}{24p_{1}}%
,\frac{1}{96}\}$ where $p_{1}$ is the index of the $L^{p_{1}}$-norm used for
the density $n(\cdot,t)$ with $t\in\lbrack\sigma,T_{max})$. It follows that
$(n,c,\zeta)$ can be extended globally in time and belongs to
(\ref{space-measure-aux-1})-(\ref{space-measure-aux-2}), for every
$T\in(0,\infty).$ By standard arguments in PDEs, the global uniqueness follows
from the local one. Thus, by using the property given in item $(ii)$ in Remark
\ref{rem-teo-measure}, we are done.\hfill\hfill$\diamond$\vskip12pt

\begin{remark}
\label{Rem-global}

\begin{itemize}
\item[$(i)$] (Persistence property) Taking in particular $p_{1}=p_{3}=1$ and
$\alpha_{1}=\alpha_{3}=0$ in (\ref{space-measure-aux-1}%
)-(\ref{space-measure-aux-2}), we obtain the solution in the persistence
space, according to Remark \ref{rem-teo-measure} $(iii)$. More precisely,
$n\in C_{w}([0,\infty);\mathcal{M}(\mathbb{R}^{2}))$, $c\in C_{w}%
([0,\infty);L^{\infty}(\mathbb{R}^{2}))$ and $\zeta\in C_{w}([0,\infty
);\mathcal{M}(\mathbb{R}^{2}))$ where the continuity in time should be meant
in the strong sense for $t>0$ and in the distributional sense at $t=0^{+}$.

\item[$(ii)$] (Regularity of solutions) Considering $\nabla\phi\in
H^{m+1}(\mathbb{R}^{2})$ with $m\in\mathbb{N}$, we obtain that the solution in
Theorem \ref{global-M} satisfy the regularity (\ref{Aux8.1})-(\ref{Aux8.2})
for all $m_{1}\in\mathbb{N}$ and $m_{2}=m+1,$ according to Remark
\ref{rem-teo-measure} $(i).$


\end{itemize}
\end{remark}


\end{document}